\documentclass[fleqn]{interact}

\usepackage{epstopdf}
\usepackage[caption=false]{subfig}
\theoremstyle{plain}

\theoremstyle{definition}

\theoremstyle{remark}

\usepackage{lineno}
\usepackage{times}

\usepackage[utf8]{inputenc}
\usepackage[T1]{fontenc}
\usepackage{bbm}
\usepackage{lmodern}
\usepackage{amssymb }
\usepackage{amsmath}
\usepackage{latexsym}
\usepackage{algorithm}	
\usepackage{algpseudocode}
\algrenewcommand\algorithmicrequire{\textbf{Input:}}
\algrenewcommand\algorithmicensure{\textbf{Output:}}
\usepackage{graphicx}
\usepackage{multicol}
\usepackage{multirow}
\usepackage{mathtools}
\usepackage{tikz}
\usepackage{tabularx}
\usepackage{letltxmacro}
\usepackage{wrapfig}
\usepackage{lipsum}
\usepackage{thmtools}
\usepackage{thm-restate}
\usetikzlibrary{positioning}
\usepackage[hyphens]{url}
\usepackage{pdfpages}
\usepackage[colorlinks]{hyperref}
\usepackage[authoryear,round]{natbib}
\makeatletter
\def\ps@pprintTitle{
	\let\@oddhead\@empty
	\let\@evenhead\@empty
	\let\@oddfoot\@empty
	\let\@evenfoot\@oddfoot
}
\def\BState{\State\hskip-\ALG@thistlm}

\begin{document}

\articletype{ARTICLE TEMPLATE}

\title{Simultaneous Planning of Liner Ship Speed Optimization, Fleet Deployment, Scheduling and Cargo Allocation with Container Transshipment}

\author{
\name{Jasashwi Mandal\textsuperscript{a}\thanks{CONTACT Jasashwi Mandal. Email: jasashwi.mandal@nitie.ac.in}, Adrijit Goswami\textsuperscript{b}, Lakshman Thakur\textsuperscript{c}, Manoj Kumar Tiwari\textsuperscript{a,d}
}
\affil{\textsuperscript{a} National Institute of Industrial Engineering, Mumbai, India; \textsuperscript{b} Department of Mathematics, Indian Institute of Technology Kharagpur, India; \textsuperscript{c} Operations and Information Management, Business School, University of Connecticut, Storrs, CT, USA;
\textsuperscript{d} Department of Industrial and Systems Engineering, Indian Institute of Technology Kharagpur, India}
}

\maketitle

\begin{abstract}
Due to a substantial growth in the world waterborne trade volumes and drastic changes in the global climate accounted for $CO_2$ emissions, the shipping companies need to escalate their operational and energy efficiency. Therefore, a multi-objective mixed-integer non-linear programming (MINLP) model is proposed in this study to simultaneously determine the optimal service schedule, number of vessels in a fleet serving each route, vessel speed between two ports of call, and flow of cargo considering transshipment operations for each pair of origin-destination. This MINLP model presents a trade-off between economic and environmental aspects considering total shipping time and overall shipping cost as the two conflicting objectives. The shipping cost comprises of $CO_2$ emission, fuel consumption and several operational costs where fuel consumption is determined using speed and load. Two efficient evolutionary algorithms: Nondominated Sorting Genetic Algorithm II (NSGA-II) and Online Clustering-based Evolutionary Algorithm (OCEA) are applied to attain the near-optimal solution of the proposed problem. Furthermore, six problem instances of different sizes are solved using these algorithms to validate the proposed model.
\end{abstract}


\section{Introduction}\label{introduction}

Liner shipping performs a significant role in the growth of world economy due to gradual increase in global operations. Around 90\% of the global trade is borne by the maritime shipping industry. United Nations Conference on Trade and Development \cite{asariotis2018review} revealed that the containerized trade is responsible for about 17.1\% of the international seaborne trade in tons. In 2018, global containerized trade increased by 2.6\% and the volume exceeded 170 million 20-foot equivalent units (TEUs). As per UNCTAD records, containerized trade has experienced an annual average growth rate of 8.1\% between 1980 and 2017. Over past two decades, the rapid development in flow containerized cargo is an outcome of rising global trade volumes, economic globalization, and more efficient container handling facilities in ports. Liner shipping is an essential marine transportation mode, where a ship owner sails the ship to visit sequentially a number of ports. Each liner shipping company runs its own services on predetermined routes each containing a set of port calls. Liner ships mainly transports containerized cargo on the regularly scheduled service routes. \par
The vessel speed in maritime shipping remains low compared to the other transportation modes. Since long voyages can take 1 to 2 months, accelerating vessel speed may result in a quick delivery of products, low storage costs and more trade throughput per unit time. On the other hand, higher sailing speed generates more $CO_2$ emission. Fuel consumption is directly proportional to vessel speed and burning more fuel will definitely increase emissions. But, whenever seaborne trades are low and fuel price is high, these profits may become difficult to achieve. Then vessels prefers to slow down, and follows slow steaming policy in such situations. Hence, slow steaming is a prevalent policy to reduce fuel costs and emissions together and also beneficial both in economic and environmental aspects. The number of vessels required on each route can be determined by fleet deployment decision. Depending on fleet size, the fixed operation cost of vessels will vary. The optimal vessel speed and the duration of port operations affect vessel arrival time at each portcall, i.e., service scheduling decision. The former variables also depend on the cargo routing decision which decides flow of containers (TEU) for all origin-destination (OD) pairs fulfilling cargo demands at each port. Container demands for an origin and destination pair belonging to different routes can be transshipped. Transshipment operations further originate handling and inventory holding cost for the containers at the corresponding ports during transshipment.  Therefore, a manager of a shipping company has to choose a feasible plan involving these interconnected decisions.
\par
Now, the strategy is to find (1) optimal speed: ship speed on each sailing leg of the routes; (2) fleet deployment decision: Number of required vessels deploying on all routes; (3) service scheduling decision: vessels' arrival time at any portcall on all routes every week; (4) cargo flow considering transshipment: Flow of containerized cargos shipped by the liner services from each origin-destination pair. In this study, we simultaneously optimize the sailing speed, service schedule, fleet deployment, and cargo flow decisions where transshipment between routes is allowed. A bi-objective mixed-integer non-linear programming (MINLP) model is proposed to solve this real and complex problem. This model captures the conflicting nature of a liner shipping problem where there are two objectives: (1) to minimize total cost which contains fixed operation cost, berth occupancy charge, transshipment handling and holding cost, loading and unloading cost, fuel consumption cost, cost of $ CO_2 $ emission while sailing in the sea and cargo handling at the ports and (2) to minimize total time consumed by the vessels to fulfill container shipment demands of each port.
\par
Innovative features of this research reside in the multiobjective model formulation that integrates various complex and crucial operational decisions of maritime industry, e.g., determining sailing speed on each voyage leg, deploying capacitated fleet, time scheduling for each port of call, and cargo flows for each origin destination pair. A major portion of the researchers addressing ship routing or scheduling problems, does not allow transshipment operations in their study although it is very convenient and economic to be followed in shipping industry. This study considers transshipment operation which allows containers to be transshipped at the ports. At the time of transshipping the containers from one vessel to another, the containers have to be stored for the time the next vessel takes to visit the transshipment port and this incurs an inventory holding cost. Another realistic factor taken into account here is ship's bunker consumption that depends on both its load and speed. Moreover, to improve environmental sustainability the model includes the costs related to carbon emissions in sea during voyage and at ports due to cargo handling operations.
\par

\section{Literature Survey}\label{Literature survey}
As suggested by \cite{meng2013containership}, there can be three categories of the decisions taken by liner shipping companies: 1. operational decisions (i.e., container routing, cargo booking, potential rejection of cargo, vessel rescheduling); 2. tactical decisions (i.e., fleet deployment, determining frequency, schedule construction, optimization of speed); and 3. strategic decisions (i.e., network design, alliance strategy, fleet size and mix). The literature survey presented herein especially focuses on a liner shipping problem with tactical and operational level decision making. This section provides the related literature review on sailing speed optimization, fleet deployment, scheduling decision, container routing and green shipping. \cite{qi2012minimizing} presented a ship scheduling problem to optimize the total bunker consumption. A robust liner ship schedule was developed by \cite{wang2012robust} with uncertain port operation time and schedule recovery. 
To stabilize the trade-off between customer service level and higher bunker consumption, \cite{brouer2013vessel} developed a vessel schedule recovery strategy. 
Later, \cite{wang2014liner} designed a tactical service scheduling problem considering port time windows. Both the scheduling problem and cargo routing problem have been addressed by some researchers.
Furthermore, a ship scheduling problem was addressed in the multi-objective MINLP model proposed by \cite{dulebenets2018comprehensive}.
Similarly, \cite{zhuge2020schedule} studied a schedule design problem faced by liner shipping companies under voluntary vessel speed reduction incentive programs (VSRIPs).
\cite{wang2021deploying} addressed the problem of jointly deploying, sequencing, and scheduling a fleet of heterogenous vessels in a shipping route.
Later, \cite{zheng2022k} developed a robust container vessel sequencing (RCVS) problem with week-dependent demands in a shipping route.
Recently, \cite{du2023container} presented a MILP problem that involves optimizing shipping schedules and selecting a shipper based on their behavior.
\par
Although there are a number of studies on fleet deployment problem, most of them assumed direct shipment without transshipping activities and known containerized cargo demand \citep{lee2012minimizing, gelareh2010novel}. \cite{meng2010chance} addressed fleet deployment problem having uncertainties in cargo shipping demand, without considering container transshipping activities with the help of a chance-constrained programming framework. A liner ship fleet deployment (LSFD) problem allowing transshipping activities was presented by \cite{meng2012liner}.
Similarly, \cite{huang2015liner} addressed a liner service network design and ship deployment problem with cargo transshipment.
An integrated shipping network was designed by \cite{zheng2015network} considering ship deployment, cargo routing and uncertain demand.
\cite{ng2018fleet} investigated a liner ship deployment problem with incomplete shipment demand information without transhipping containers. Moreover, \cite{zhen2019fleet} developed a demand fulfillment problem considering overload risk of containers, port resources and vessel size while transshipment is allowed. Similarly, \cite{zhen2019integrated} again studied a ship deployment problem along with cargo allocation and ship scheduling considering transshipment operations. 
Furthermore, a two-stage robust optimization model has been suggested by \cite{lai2022robust} for managing the deployment of ship fleets and revenue in a liner shipping network while taking into account demand uncertainties.
\par
Several studies on liner shipping considered that vessels travel at a predetermined speed \citep{gelareh2010novel,gelareh2011fleet,meng2012liner}. \cite{ronen2011effect} highlighted the effect of sailing speed on operating cost and addressed the problem of speed optimization along with the number of vessels required and service frequency. 
At the strategic level, \cite{xia2015joint} addressed speed optimization, cargo allocation, and fleet deployment jointly to maximize total gain of shipping companies. They first considered a bunker consumption function that depends on load and speed, and then approximated the function to make the model linear. But they did not consider transshipments among different service routes and also the enviromental aspect is not captured in their study.
Later, \cite{karsten2018simultaneous} developed a decision support tool for profit maximization of shipping companies finding optimal number of vessels required and vessel speed for a global network. To optimize the service level and total fuel consumption in a liner shipping problem with stochastic time windows and port times, a dynamic programming formulation finding suitable sailing speed was developed by \cite{aydin2017speed}.
Moreover, \cite{koza2019liner} proposed a branch-and-price search technique to address a liner service scheduling and container routing problem incorporating bunker consumption that is a function of load and speed.
\cite{zhen2020route} proposed a bi-objective mixed integer linear programming model, aiming to optimize sailing routes and speeds within and outside the ECA while minimizing the total fuel cost and emissions.
Thereafter, \cite{brahimi2021exact} presented an exact algorithm to solve a single liner service design problem with speed optimization. In a recent study, \cite{wang2022probabilistic} addressed a liner shipping problem that involved optimizing vessel speeds to minimize the overall costs of the fleet. These costs included operating expenses, capital expenditures, and voyage-related expenses.
\par
Very few studies have jointly considered major decisions of maritime shipping industry capturing the environmental issues. Due to the considerable growth in the international seaborne trade volume, $ CO_2 $ emission has become an alarming issue; despite of this fact environmental concern has not attracted much attention in the literature except few researchers. For example, \cite{wen2017multiple} aimed at designing a multiple ship speed optimizing and routing model under cost, time and environmental objectives. \cite{fagerholt2015two} optimized speed of ships sailing in and out of Emission Control Areas (ECAs) that have restriction on sulfur emissions.
A green vessel scheduling problem is presented by \cite{dulebenets2018greenijtst} considering $ CO_2 $ emission costs at each port of call and also during voyages in sea.
Another interesting study by \cite{wang2021liner}, addressed a liner shipping service planning problem that integrates fleet deployment, schedule design, and sailing path and speed optimization, considering the effect of ECAs. But they have not considered container transshipment in the liner service.
\par
Therefore, to the best of our knowledge, for the first time in literature this present study aims to address a holistic liner shipping problem and jointly take the complex decisions of speed optimization, scheduling, fleet deployment and cargo flow with container transshipment under cost (including environmental aspect) and time objectives. Hence, the major contribution of the study is trifold. Firstly, it finds a trade-off between the two conflicting objectives: 1. total cost incurred during shipping service where the environmental impact is also presented as a cost, and 2. total journey time of vessels deployed to services. The second contribution is considering container transshipment. Generally, a liner shipping company may fail to provide direct services between any two ports since large number of ports are involved in global shipping services. Many liner shipping problems does not consider container transshipment operations to avoid complexity in the model. Third, the bunker consumption here is taken as a function of both speed and load, where a limited number of studies has taken load as a factor for determining fuel consumption. In summary, this study has proposed a multi-objective model in maritime transportation focusing on sailing speed, fleet deployment, cargo flow, service scheduling, and environmental impact of ships.

\section{Problem description and Mathematical formulation}\label{Problem description and Mathematical formulation}
Let us consider that a number of liner service routes (represented by $ R $, the set of routes) covering a number of ports (set of ports is denoted by $ P $) are operated by a liner shipping company. Since liner shipping services in general follows weekly pattern, the service frequency of each ship route $ r \in R $ is maintained weekly. Let $ I_r $ be the set of ports of call (or set of legs) on vessel route $ r $ i.e. vessel route $ r $ has $ | I_r| $ legs as well as $ |I_r| $ ports of call. Any liner service route $ r \in R $ is explained as a port rotation: $ p_{r1} \rightarrow p_{r2} \rightarrow ...\rightarrow p_{r|I_r|} \rightarrow p_{ri} $, where $ p_{ri} $ is the $ i^\mathrm{th} $ port of call and $ i={1,2,...,|I_r|} $. The sequence of ports to be maintained (port rotation) is given. The voyage between port $ p_{ri} $ and port $ p_{r,i+1} $ is said to be leg $ i $ , with $ p_{r,|I_r|+1} = p_{r1}$. Usually, the shipping company determines the port sequence at the strategic level. Each route deploys a fleet of vessels to maintain weekly service. Generally, ships deployed by the shipping company are either their own ships or chartered ships. These ships have different categories based on load capacities, fixed operating costs, berth occupancy charges, ship-weight and other characteristics related to ship. Moreover, the sailing time for each leg can be set by controling the speed and is also related to the decision of fleet deployment on the routes. The vessel arrival time at each port is influenced by the vessel speed on each leg. $\theta_{rir'i'}$ denotes the time difference (in days) between vessel arrival times at port of call ($ r,i $) and port of call ($ r',i' $), i.e., the time required to store containers at the transsipment port. A transshipment from the $ (r,i) $ portcall to the $ (r',i') $ portcall is described as ($ r,i,r',i' $), that indicates that for route $ r $ the $ i^\mathrm{th} $ portcall and for route $r'$ the $ {i'}^\mathrm{th} $ portcall is associated to the same shipping port in the service, i.e., $ p_{ri} = p_{r'i'} $. The multiobjective liner shipping problem proposed here takes the following into consideration:\\
\begin{enumerate}
	\item The port rotation on each route in the shipping network is provided.
	\item For each origin-destination, volume of container shipment demand (TEU per week) is fixed.
	\item The deployed vessels are heterogeneous in terms of cost structure and capacity.
	\item Each vessel will start from a known port i.e., its initial position.
	\item Bunker consumption on each voyage leg is based on both the vessel speed and payload.
	\item To estimate the total amount of carbon emission, Marine diesel oil (MDO) and Heavy fuel oil (HFO) carbon emission coefficients are considered.
\end{enumerate}
\subsection{Model Formulation}\label{class}
\subsubsection{Sets and Indices}
The sets and indices used here are enlisted below:
\begin{table}[H]
	\begin{tabularx}{\textwidth}{lX}
		$ R \;\;\;\;\;\;\;\;\;\;\;\;$ & Set of shipping routes;\\
		$ r \;\;\;\;\;\;\;\;\;\;\;\;\;\;\;\;$ & Index of each shipping route, $r\in R$;\\
    $ i $ & Index of each leg $ (i,i + 1) $ or each port of call on a shipping route;
  	\end{tabularx}
	\end{table}
	\begin{table}[H]
		\begin{tabularx}{\textwidth}{lX}
    $ I_r $ & Set of all legs (or port of calls) on shipping route $ r $ ;\\
		$ P $ & Set of all ports;\\
		$ p $ & Index of each port, $ p\in P $;\\
		$p_{ri}$ & Port index corresponding to the $ i^\mathrm{th} $ port of call on shipping route $ r $, $p_{ri}\in P$;\\
		$(r,i,r',i')\;\;\;\;$ & Transshipment index from port $p_{ri}$ to port $p_{r'i'}$;
		$ p_{ri} = p_{r'i'} $ where $r,r'\in R;i\in I_r,i'\in I_{r'}$;  ;\\
		$ ST $ &Set of transshipment quadruples $ (r, i, r', i') $; $ ST=\{(r,i,r',i')|\;p_{ri}=p_{r'i'}\}$;\\
		$I_{rp}$ & Set of legs (or ports of call) of port $p$ on liner service route $ r $ ;\\
		$ R_p $ & Set of all shipping routes covering port $p$;\\
		$ V $ & Set of different types of vessels;\\
		$ v $ & Index of vessel type, $ v\in V $;\\
		TEU & Twenty-foot equivalent unit (unit of cargo capacity).
	\end{tabularx}
\end{table}
\vspace{-0.6cm}
\subsubsection{Parameters}
The parameters for the proposed problem are enlisted below:
\vspace{-0.2cm}
\begin{table}[H]
	\begin{tabularx}{\textwidth}{lX}
		$ C^{opr}_v $ & Operating cost of each vessel type $ v\in V $ (USD/week). It includes the costs of repair and maintenance, crew, stores, lubricants, and insurance.
		It does not depend on number of trips and the cost will incur until the ship stops to provide service.\\
		$  C^{berth}_{v} $ & Charges to occupy berth (USD/hour) for each vessel type $ v\in V$. It depends on the berth occupancy time and rate of vessel charges for each type at each port;\\
		$  C^{fix}_{rv} $ & Costs related to voyages on route $r$ for vessel type $v$ (USD/week). It includes fixed costs for calling at ports and canal dues;\\
		$ C^{load}_p $ & Charges for loading or charging containers to ships at port $ p\in P $ (USD/TEU);\\
		$ C^{disc}_p $ & Charges for unloading or discharging containers from ships at port $ p\in P $ (USD/TEU);\\
		$ C^{trans}_p $ & Cost for handling transshiped containers at port $ p\in P $ (USD/TEU).\\
		$ C^{hold}_p $ & Holding cost of containers transshipped at port $ p\in P $ (USD/TEU/hour);\\
		$ d_{od} $ & Container shipment demand need to be transported from one port $o\in P$ to other port $ d\in P$ in a week (TEUs/week);\\
		$Cap_v$ & Capacity of a vessel type v;\\
		$ T_{pv}\;\;\;\;\;\;\;\;\;\;\;\;\;\;$ &  Average handling time of one TEU container at port $ p $ for vessel type $ v $ (hour/TEU);\\
		$ L_{ri} $ & Voyage length of $ i^\mathrm{th} $ leg of shipping route $ r $ (nautical miles);\\
		$n^{min}_r$, $n^{max}_r$ & Minimum and maximum number of vessels which can be employed on shipping route $ r $ respectively;\\	
		$u^{min}_{ri}$, $u^{max}_{ri}$ & Minimum and maximum speed of vessels on $ i^\mathrm{th} $ leg on shipping route $ r $ respectively.\\
		$ k_{prir'i'} $ & It has the value 1 if and only if port $ p $ does transshipment $(r, i, r',i') $, and otherwise 0;\\
  	$  F_{riv} $ & Fuel consumption on $ i^\mathrm{th} $ leg of shipping route $r$ for vessel type $v$;\\
  $c^{F}$ & Cost of unit fuel consumption (USD/ton);\\
  		$c^{E}$ & Cost of $ CO_2 $ emission;\\
    		$E_{sea}$ & $ CO_2 $ emission factor at sea (tons of $ CO_2 $/ton of fuel);\\
      			$E_{port}\;\;\;\;\;\;\;\;\;\;$ & $ CO_2 $ emission factor at ports (tons of $ CO_2 $/TEU handled).
   	\end{tabularx}
		\vspace*{-1cm}
	\end{table}
	\vspace*{-0.4cm}
\subsubsection{Variables}
The decision variables are enlisted as follows:
\vspace{-0.2cm}
\begin{table}[H]
	\begin{tabularx}{\textwidth}{lX}
		$ n_r $ & Number of ships required to be deployed on shipping route $ r\in R $;\\
		$ x_{rv} $ & 1, \hspace{0.3cm} if shipping route $ r $ is using vessel type $ v\in V $;\\
		& 0, \hspace{0.3cm} otherwise;\\
		$t_{ri}$ & Time (hours) when a vessel arrives at the $i^\mathrm{th}$ portcall on shipping route $ r $ with $ i = 1,...,|I_r|+1 $; Let us take $ 0 \le t_{r1} \le 144$ (without loss of generality); $ t_{r,|I_r|+1} $ denotes the time when the vessel comes back to the $1^\mathrm{st}$ portcall on shipping route $ r $ , i.e., $ t_{r,|I_r|+1} $ equals $t_{r1}$ added to the time (in hours) needed to finish a round-trip voyage by vessels;\\
		$\theta_{rir'i'}\;\;\;\;\;\;\;\;\;\;$ & Difference between vessel arrival times (hours) at portcall $ (r,i) $ and portcall $ (r',i') $;\\
		$u_{ri}$ & Sailing speed (knots) on $i^\mathrm{th}$ leg of shipping route $r$;\\
		$z^{load}_{rip}$ & Number of containers starting from port $ p\in P $ and loaded or charged at  $i^\mathrm{th}$ portcall on liner service route $ r\in R $ (TEUs/week);\\
		$z^{disc}_{rip}$ & Number of containers starting from port $ p\in P $ and unloaded or discharged at $i^\mathrm{th}$ portcall on liner service route $ r\in R$ (TEUs/week);\\
		$f_{rip}$ & Number of containers starting from port $ p\in P $ and charged in the vessels that sails on $i^\mathrm{th}$ leg of shipping route $ r\in R $ (TEUs/week);\\
		$\gamma_{rir'i'}$ & An integer associated to variable $\theta_{rir'i'}$ and used to change the difference between arrival times $t_{ri}$ and $t_{r'i'}$ to a non-negative integer.
	\end{tabularx}
	\vspace*{-0.2cm}
\end{table}	
\subsubsection{Mathematical model}
To address the proposed problem, a multi-objective mathematical model is formulated. The objective of the model is to determine optimal sailing speed, fleet depolyment, vessel  schedule and flow of cargo fulfilling the demand so that total shipping cost and shipping time are minimized. Eq.\eqref{eq1} outlines the first objective function that depicts the total cost related to different shipping operations and complete the round-trip. The first objective function consists of eight components. First term is the cost related to voyages and operating vessels, and berth occupancy charge for vessels is presented by the second component. Next, third one is associated with handling cost for transshipped containers. Fourth term depicts the transshipment holding cost depending on flow of containers and ship schedule. Total discharge and loading costs are described by the fifth component, and the sixth term presents total weekly fuel consumption cost incurred while sailing in sea. Next, the seventh and eighth part provide total $ CO_2 $ emission cost associated with the fleet of ships while flowing in sea and serving at ports respectively.
Eq.\eqref{eq2} formulates the second objective that depicts the total time spent to complete the round-trip journey.
These two objectives conflicting in nature, are described below:\vspace{0.1cm}\\
\emph{First objective}:\\
Minimize total cost = Ship and voyage operating cost + Berth occupancy charge + Transshipped container handling cost + Transshipment holding cost + Loading and unloading cost + Fuel consumption cost + $ CO_2 $ emission cost while sailing in the sea + $ CO_2 $ emission cost while operating at the port.\vspace{0.1cm}\\
\emph{Components of the first objective}:\\
\vspace{-0.4cm}
\begin{flalign}
	1. \text{ Ship and voyage operating cost = }\sum_{r\in R}\sum_{v\in V}\big(n_rC^{opr}_v + C^{fix}_{rv}\big)x_{rv}\notag
\end{flalign}
\vspace{-0.6cm}
\begin{align}
	2. \text{ Berth occupancy charge = }\sum_{r\in R}\sum_{i\in I_r}\sum_{o\in P}\sum_{v\in V} C^{berth}_{v}T_{p_{ri}v}x_{rv}\big(z^{load}_{rio}+z^{disc}_{rio}\big)\notag
\end{align}
\vspace{-0.6cm}
\begin{align}
	&3. \text{ Transshipped container handling cost = } \notag\\
 &\frac{1}{2}\sum_{p\in P}C^{trans}_p\bigg(\sum_{r\in R_p}\sum_{i\in I_{rp}}\sum_{\substack{o\in P\\o \ne p}}\big(z^{load}_{rio}+z^{disc}_{rio}\big)-\sum_{d\in P}d_{pd}-\sum_{o\in P}d_{op}\bigg)\notag
\end{align}
\vspace{-0.6cm}
\begin{align}
	&4. \text{ Transshipment holding cost = }\notag\\
 &\frac{1}{2}\sum_{p\in P}C^{hold}_p \hspace{-0.3cm} \sum_{(r,i,r',i')\in ST}\hspace{-0.5cm}k_{prir'i'}\theta_{rir'i'}\bigg(\sum_{r\in R_p}\sum_{i\in I_{rp}}\sum_{\substack{o\in P\\o \ne p}}\big(z^{load}_{rio}+z^{disc}_{rio}\big)
	-\sum\limits_{d\in P}d_{pd}-\sum_{o\in P}d_{op}\bigg)\notag
\end{align}
\vspace{-0.6cm}
\begin{align}
	5. \text{ Loading and unloading cost = }\sum_{o\in P}\sum\limits_{d\in P}\big(C^{load}_o+C^{disc}_d\big)d_{od}\notag
\end{align}
\vspace{-0.4cm}
\begin{align}
	6. \text{ Fuel consumption cost = }c^F\sum_{r\in R}\sum_{i\in I_r}\sum_{v\in V}\frac{L_{ri}}{24u_{ri}}F_{riv}x_{rv}\notag
\end{align}
\vspace{-0.3cm}
\begin{align}
	7. \text{ $ CO_2 $ emission cost while sailing in the sea = }c^E E_{sea}\sum_{r\in R}\sum_{i\in I_r}\sum_{v\in V}\frac{L_{ri}}{24u_{ri}}F_{riv}x_{rv}\notag
\end{align}
\vspace{-0.4cm}
\begin{align}
	& 8. \text{ $ CO_2 $ emission cost while operating at the port = } \notag\\
 & c^E E_{port}\sum_{r\in R}\sum_{i\in I_r}\sum_{o\in P}\sum_{v\in V}x_{rv}\big(z^{load}_{rio}+z^{disc}_{rio}\big)\notag
\end{align}
\emph{Second objective}:
\vspace{-0.2cm}
\begin{align}
	\text{Minimize total time =}\sum_{r\in R}\sum_{v\in V}x_{rv}\sum_{i\in I_{r}}\bigg(\frac{L_{ri}}{u_{ri}}+T_{p_{ri}v}\sum_{o\in P}\big(z^{load}_{rio}+z^{disc}_{rio}\big)\bigg)\notag
\end{align}
The multi-objective problem mentioned above is framed as a MINLP model provided below:
\begin{align}\label{eq1}
	&\text{Minimize }F_1=
	\sum_{r\in R}\sum_{v\in V}\big(n_rC^{opr}_v + C^{fix}_{rv}\big)x_{rv} + \sum_{r\in R}\sum_{i\in I_r}\sum_{v\in V} C^{berth}_{v}T_{{p_{ri}v}}x_{rv}\sum_{o\in P}\big(z^{load}_{rio}+z^{disc}_{rio}\big) \notag \\
	&+\frac{1}{2}\sum_{p\in P}C^{trans}_p\bigg(\sum_{r\in R_p}\sum_{i\in I_{rp}}\sum_{\substack{o\in P\\o \ne p}}\big(z^{load}_{rio}+z^{disc}_{rio}\big)-\sum_{d\in P}d_{pd}-\sum_{o\in P}d_{op}\bigg)\notag\\
	&+\frac{1}{2}\sum_{p\in P}C^{hold}_p\sum_{(r,i,r',i')\in T}k_{prir'i'}\theta_{rir'i'}\bigg(\sum_{r\in R_p}\sum_{i\in I_{rp}}\sum_{\substack{o\in P\\o \ne p}}\big(z^{load}_{rio}+z^{disc}_{rio}\big)-\sum\limits_{d\in P}d_{pd}-\sum_{o\in P}d_{op}\bigg)\notag\\
	&+\sum_{o\in P}\sum\limits_{d\in P}\big(C^{load}_o+C^{disc}_d\big)d_{od} + c^F\sum_{r\in R}\sum_{i\in I_r}\sum_{v\in V}\frac{L_{ri}}{24u_{ri}}F_{riv}x_{rv} \notag\\
	& +c^E E_{sea}\sum_{r\in R}\sum_{v\in V}\sum_{i\in I_r}\frac{L_{ri}}{24u_{ri}}F_{riv}x_{rv}+ c^E E_{port}\sum_{r\in R}\sum_{i\in I_r}\sum_{o\in P}\sum_{v\in V}x_{rv}\big(z^{load}_{rio}+z^{disc}_{rio}\big)
\end{align}
\begin{align}\label{eq2}
	\text{Minimize }F_2=
	\sum_{r\in R}\sum_{v\in V}x_{rv}\sum_{i\in I_{r}}\bigg(\frac{L_{ri}}{u_{ri}}+T_{p_{ri}v}\sum_{o\in P}\big(z^{load}_{rio}+z^{disc}_{rio}\big)\bigg)
\end{align}
Fuel consumption function based on both payload and speed can be practically estimated by the  modified admiralty formula \citep{janic2014advanced, barrass2004ship, psaraftis2014ship}: $F_{riv}=k_v.u_{ri}^3.(\Delta_v)^{2/3} $, where $\Delta_v$ is the weight of the water the vessel of type $ v $ displaces (denoted as the displacement of a vessel type $ v $). $ k_v>0 $ are vessel dependent parameters. Vessel displaces the weight equal to sum of the vessels' actual deadweight (dwt) and lightweight (lwt). The deadweight defines the exact weight of everything the ship is carrying i.e., the total weight of fresh water, ballast water, fuel, loaded containerized cargo, provisions and crew. The weight of an fully empty ship is said to be its lightweight. Here, the definition is considered as $ \Delta_v := lwt_v + owt_v + w_v $ for a vessel type $ v\in V $, where $ lwt_v $ is a vessel's lightweight and $ owt_v $ represents the average weight of fresh water, ballast water, fuel, provisions and crew. $ w_v $ is the weight of all loaded containers. Since the change in weight of fresh water, ballast water and fuel is in general small unlike a vessel's total deadweight, $ owt_v $ is assumed to be constant for simplification. Hence, as suggested by \cite{koza2019liner}, total fuel consumption ($ F_{riv} $) can be defined as a payload and speed dependent function, $F_{riv}=k_v.u_{ri}^3.(lwt_v + owt_v + w_v)^{2/3} $ for a vessel type $ v\in V $.
\vspace{0.1cm}\\
\emph{Constraints}:
\begin{align}
	& n^{min}_r \le n_r \le n^{max}_r,\;\;\;\;\forall r\in R \label{eq3}\\
	& \sum_{v\in V} x_{rv}=1,\;\;\;\;\forall r\in R \label{eq4}\\
	& u^{min}_{ri} \le u_{ri}\le u^{max}_{ri},\;\;\;\;\forall r\in R,\; \forall i\in I_r \hspace{8cm} \label{eq5}
\end{align}
Constraints \eqref{eq3} restrict the total number of vessels required to operate on each sailing route taking minimum and maximum vessels as $ n^{min}_r $ and $ n^{max}_r $ respectively. Constraints \eqref{eq4} necessitate the condition that each shipping route can deploy exactly one vessel class. Constraints \eqref{eq5} depicts the span of vessel speed on each voyage leg of the shipping routes considering minimum and maximum speed level as $ u^{min}_{ri} $ and $ u^{max}_{ri} $ respectively.
	\begin{align}
		& \sum_{o\in P}f_{rio}-\sum_{v\in V}Cap_v x_{rv} \le 0, \;\;\;\;\forall r\in R,\; \forall i\in I_r \label{eq6}\\
		& \sum_{r\in R_d}\sum_{i\in I_{rd}}\big(z^{disc}_{rio}-z^{load}_{rio}\big)=d_{od},\;\;\;\;\forall d\in P,\;\forall o\in P,\;d\ne o \hspace{5cm}\label{eq7}\\
		& f_{r,i-1,o}+z^{load}_{r,i,o}=f_{rio}+z^{disc}_{r,i,o}, \;\;\;\;\forall r\in R,\; \forall i\in I_r,\;\forall o\in P \label{eq8}
	\end{align}
Constraints \eqref{eq6} enforce that the cargo flow on each voyage leg of the shipping routes must be within the selected vessel capacity. Constraints \eqref{eq7} imply that the cargo shipping demand is fulfilled for all ports and the cargo flow on all voyage legs and all portcalls on every liner service route is represented by constraints \eqref{eq8}.
\begin{align}
	& 0 \le t_{r1} \le 144,\;\;\;\;\forall r\in R\label{eq9}\\
	& t_{r,i+1}=t_{ri}+\sum_{v\in V}\bigg(T_{p_{ri}v}\sum_{o\in P}\big(z^{load}_{rio}+z^{disc}_{rio}\big) + \frac{L_{ri}}{u_{ri}}\bigg)x_{rv},\;\;\;\;\forall r\in R,\; \forall i\in I_r \label{eq10}
\end{align}
Constraints \eqref{eq9} impose that within seven days of the period of planning the first vessel employed on each shipping route serves the $ 1^{st} $ port of call on that route and it leads the liner shipping services to follow the weekly frequency. Next, the vessel visiting time at any portcall and at its subsequent portcall in a liner service route are connected by constraints \eqref{eq10}.
\begin{align}
	& t_{r'i'}-t_{ri} + 168\gamma_{rir'i'}=\theta{rir'i'},\;\;\;\;\forall (r,i,r',i')\in ST \label{eq12}\\
	& 0 \le \theta_{rir'i'} \le 144, \;\;\;\;\forall (r,i,r',i')\in ST \label{eq13}\\
	& -n_{r'} \le \gamma_{rir'i'} \le n_r, \;\;\;\;\forall (r,i,r',i')\in ST \label{eq14}
\end{align}
The difference between vessel arrival times at portcall $ (r,i) $ and $ (r',i') $ involved in transshipment $ (r,i,r',i') $ is determined by constraints \eqref{eq12}-\eqref{eq14}.
\begin{align}
	& 168n_r \ge \sum_{v\in V}x_{rv}\sum_{i\in I_{r}}\Big(\frac{L_{ri}}{u_{ri}}+T_{p_{ri}v}\sum_{o\in P}\big(z^{load}_{rio}+z^{disc}_{rio}\big)\Big), \;\;\;\;\forall r\in R \label{eq15}\\
	&f_{rio}=0, \;\;\;\; \forall r\in R,\; \forall i\in I_r,\;o=p_{r,i+1} \label{eq16}\\
	&z^{disc}_{rio}=0, \;\;\;\; \forall r\in R,\; \forall i\in I_r,\;o=p_{ri} \label{eq17}
\end{align}
Constraints \eqref{eq15} enforce that the total number of vessels serving each shipping route is sufficient to sustain a weekly pattern. Constraints \eqref{eq16} ensure that containerized cargos should not return to their origin port $ o$. Constraints \eqref{eq17} ensure that containerized cargos should not be discharged at their origin port $ o $.
\begin{align}
	&x_{rv}\in \{0,1\}, \;\;\;\;\forall r\in R, \;\forall v\in V \label{eq18}\\
	&n_r\in \mathbb{Z}^+, \;\;\;\; \forall r\in R\label{eq19}\\
	&z^{disc}_{rio}\ge 0, \;\;\;\; \forall r\in R,\; \forall i\in I_r,\; o\in P\label{eq20}\\
	&z^{load}_{rio}\ge 0, \;\;\;\; \forall r\in R,\; \forall i\in I_r,\; o\in P \label{eq21}\\
	&f_{rio}\ge 0, \;\;\;\; \forall r\in R,\; \forall i\in I_r,\; o\in P \label{eq22}\\
	&\gamma_{rir'i'}\in \mathbb{Z},\;\theta_{rir'i'}\in \mathbb{R}^+, \;\;\;\; \forall (r,i,r',i')\in ST \label{eq23}\\
	& t_{ri}, \; u_{ri}\in \mathbb{R}^+, \;\;\;\; \forall r\in R,\; \forall i\in I_r \label{eq24}
\end{align}
Lastly, constraints \eqref{eq18} - \eqref{eq24} describe the integer attributes and non-negativity of the decision variables. The above-discussed multi-objective model incorporates non-linear overall cost and time function as first and second objective respectively. It designs various decision variables including continuous, integer and binary as well as real-life constraints (both linear and non-linear) like vessel capacity, cargo flow to fulfill demand, service scheduling and required vessels.
\subsection{Linearization}
The proposed mathematical model in the previous section is a MIP model with some non-linear expressions such as $ n_rx_{rv} $, $ x_{rv}/u_{ri} $, $x_{rv}(z^{load}_{rio}+z^{disc}_{rio})$, $\theta_{risj}(z^{load}_{rio}+z^{disc}_{rio})$ and $F_{riv}x_{rv}/(24u_{ri}) $. We can transform it to a MILP model that can be solved in a commercial software based on conventional techniques.
To linearize the previous formulation, some new sets, variables and constraints are added which are given below:
\subsubsection{Newly defined sets and indices}
\vspace{-0.4cm}
\begin{table}[H]
	\begin{tabularx}{\textwidth}{lX}
		$A_{ri}\;\;\;\;\;\;\;\;\;\;\;\;$ & Set of vessel speeds on $ i^\mathrm{th} $ leg on shipping route $ r $ \\
		$ \alpha $ &  Index of vessel speeds; $ \alpha \in A_{ri} $\\
		$ H $ & Set of hours in a week; $ H=\{1,2,...,168\} $\\
		$ h $ & Index of hours in a week
	\end{tabularx}
\end{table}
\vspace{-0.3cm}
\subsubsection{Newly defined variables}
\vspace{-0.3cm}
\begin{table}[H]
	\begin{tabularx}{\textwidth}{lX}
		$n_{rv}$ & Number of $ v $ type vessels deployed on shipping route $ r $\\
		$ z_{riv} $ & Total containers loaded and unloaded at $ i^\mathrm{th} $ portcall on shipping route $ r $ for vessel type $ v $ \\
		$ \psi_{ri\alpha} $ & A binary variable, that equals one if and only if $ u_{ri} = \alpha $ for leg $ i $ on shipping route $ r $\\
		$\phi_{riv\alpha} $ & A binary variable to replace $ \psi_{ri\alpha}x_{rv} $ \\
		$\lambda_{rir'i'h} $ & A binary variable, that equals one if and only if $ \theta_{rir'i'}= h $\\
		$\delta_{prir'i'h} \;\;\;\;\;$ & Equals to the product $\displaystyle  \theta_{rir'i'}\sum_{r\in R_p}\sum_{i\in I_{rp}}\sum_{\substack{o\in P\\o \ne p}}(z^{load}_{rio}+z^{disc}_{rio}) $ if and only if $ k_{prir'i'}=1 $ and $ \theta_{rir'i'}=h $, and 0 otherwise
	\end{tabularx}
\end{table}	
\vspace{-0.4cm}
\subsubsection{Linearized objective function}	
\begin{align}\label{eq25}
	&\text{Minimize }F_1=
	\sum_{r\in R}\sum_{v\in V}C^{opr}_v n_{rv} + \sum_{r\in R}\sum_{v\in V}C^{fix}_{rv}x_{rv} + \sum_{r\in R}\sum_{i\in I_r}\sum_{v\in V} C^{berth}_{v}T_{{p_{ri}v}}z_{riv}\notag\\
	&+\frac{1}{2}\sum_{p\in P}C^{trans}_p\Big(\sum_{r\in R_p}\sum_{i\in I_{rp}}\sum_{\substack{o\in P\\o \ne p}}\big(z^{load}_{rio}+z^{disc}_{rio}\big)-\sum_{d\in P}d_{pd}-\sum_{o\in P}d_{op}\Big)\notag\\
	&+\sum_{p\in P}C^{hold}_p\sum_{(r,i,r',i')\in ST}\sum_{h\in H}\delta_{prir'i'h}+\sum_{o\in P}\sum\limits_{d\in P}\big(C^{load}_o+C^{disc}_d\big)d_{od}\notag\\
	&+\big(c^F+ c^{E}E_{sea}\big)\sum_{r\in R}\sum_{i\in I_r}\sum_{v\in V}\sum_{\alpha\in A_{ri}}\frac{L_{ri}}{24}k_v\alpha^2 w_{riv\alpha}+ c^{E}E_{port}\sum_{r\in R}\sum_{i\in I_r}\sum_{v\in V}z_{riv}
\end{align}
\vspace{-0.5cm}
\begin{align}\label{eq26}
	\text{Minimize }F_2=
	\sum_{r\in R}\sum_{v\in V}\sum_{i\in I_{r}}\sum_{\alpha\in A_{ri}}L_{ri} \frac{\phi_{riv\alpha}}{\alpha}+\sum_{r\in R}\sum_{v\in V}\sum_{i\in I_{r}}T_{p_{ri}v}z_{riv}
\end{align}
After linearization, the two contradictory objectives of minimizing cost incurred and shipping time are represented formulating the equations \eqref{eq25} and \eqref{eq26} respectively. In objective \eqref{eq1}, $ n_{r}x_{rv}C^{opr}_v $ in the first term is linearized to $ n_{rv}C^{opr}_v $ in the new objective function \eqref{eq25} using constraints \eqref{eq29}. The next non-linear term $ C^{berth}_{v}T_{p_{ri}v}x_{rv}\sum_{o\in P}(z^{load}_{rio}+z^{disc}_{rio}) $ in objective \eqref{eq1} is linearized to $ C^{berth}_{v}T_{p_{ri}v}z_{riv} $ in objective function \eqref{eq25} using constraints \eqref{eq45} and \eqref{eq46}. Now, the fourth term representing the transshipment holding cost in objective \eqref{eq1} is replaced by the linear expression $ \sum_{(r,i,r',i')\in ST}\sum_{h\in H}\delta_{prir'i'h} $ to form objective \eqref{eq25}. 	In order to linearize the model, the vessel speed is discretized into $ \{\alpha: \alpha \in A_{ri}\} $ and the set of speed levels is $ A_{ri} $. The fuel consumption function is taken as: $F_{riv}=k_v.u_{ri}^3.(lwt_v + owt_v + w_v)^{2/3} $. Hence, the expression $L_{ri}F_{riv}x_{rv}/(24u_{ri}) $ in the sixth and seventh term of objective \eqref{eq1} becomes $ L_{ri}k_v.u_{ri}^2.(lwt_v + owt_v + w_v)^{2/3}x_{rv}/24 $, which is then converted to $L_{ri}k_v\alpha^2\phi_{riv\alpha}(lwt_v+owt_v+w_v)^{2/3}/24 $ using constraints \eqref{eq40}-\eqref{eq44}.
\cite{koza2019liner} proposed a linear approximation approach to linearize the nonlinear expression $ (lwt_v + owt_v + w_v)^{2/3} $ into $ a_v + b_v. w_v $, with $ a_v $ and $ b_v $ are vessel dependent parameters. As the resulting term still $L_{ri}k_v\alpha^2\phi_{riv\alpha}(a_v + b_v. w_v)/24 $ remains nonlinear, another variable $ w_{riv\alpha} $ is introduced satisfying constraints \eqref{eq61}.
Again, the eighth term is similarly linearized to $ c^{E}E_{port}\sum_{r\in R}\sum_{i\in I_r}\sum_{v\in V}z_{riv} $ as the $ CO_2 $ emission cost while the ships being operated at ports. In the objective function \eqref{eq2}, the non-linear terms are same as discussed above and hence, they are converted to linear proceeding in the same way to form objective function \eqref{eq26}.
\subsubsection{Linearized constraints}
\vspace{-0.5cm}
\begin{align}
	& \sum_{v\in V} x_{rv}=1,\;\;\;\;\forall r\in R \label{eq27}\\
	& n^{min}_r \le n_{rv} \le n^{max}_r,\;\;\;\;\forall r\in R,\;\forall v\in V, \label{eq28}\\
	&n_{rv} \le M_1x_{rv},\;\;\;\;\forall r\in R,\;\forall v\in V, \label{eq29}\\
	& u^{min}_{ri} \le u_{ri}\le u^{max}_{ri},\;\;\;\;\forall r\in R,\; \forall i\in I_r \label{eq30}	
\end{align}
Constraints \eqref{eq27} necessitate the condition that each shipping route can deploy exactly one vessel class. Constraints \eqref{eq28}-\eqref{eq29} restrict the total number of vessels required to operate on every liner service route and define the variable $n_{rv}$. Constraints \eqref{eq30} limit the range of the vessel speed on each leg of routes considering minimum and maximum speed levels.	
\begin{align}
	& \sum_{o\in P}f_{rio}-\sum_{v\in V}Cap_v x_{rv} \le 0, \;\;\;\;\forall r\in R,\; \forall i\in I_r \label{eq31}\\
	& \sum_{r\in R_d}\sum_{i\in I_{rd}}\big(z^{disc}_{rio}-z^{load}_{rio}\big)=d_{od},\;\;\;\;\forall o\in P,\;\forall d\in P,\;o\ne d \hspace{5cm} \label{eq32}\\
	& f_{r,i-1,o}+z^{load}_{rio}=f_{rio}+z^{disc}_{rio}, \;\;\;\;\forall r\in R,\; \forall i\in I_r,\;\forall o\in P \label{eq33}
\end{align}
Constraints \eqref{eq31} enforce that the cargo flow on each voyage leg of the shipping routes must be within the selected vessel capacity. Constraints \eqref{eq32} imply that the cargo shipping demand is fulfilled for all ports. The cargo flow on all voyage legs and all portcalls on every liner service route is represented by constraints \eqref{eq33}.
\begin{align}
	& 0 \le t_{r1} \le 144,\;\;\;\;\forall r\in R\label{eq34}\\
	& t_{r,i+1}=t_{ri}+\sum_{v\in V}\Big(T_{p_{ri}v}z_{riv} + \sum_{\alpha\in A_{ri}}L_{ri} \frac{\phi_{riv\alpha}}{\alpha}\Big),\;\;\;\;\forall r\in R,\; \forall i\in I_r \label{eq35}\\
	& t_{r'i'}-t_{ri} + 168\gamma_{rir'i'}=\theta_{rir'i'},\;\;\;\;\forall\; (r,i,r',i')\in ST \label{eq37}\\
	& 0 \le \theta_{rir'i'} \le 144, \;\;\;\;\forall\; (r,i,r',i')\in ST \label{eq38}\\
	& -\sum_{v\in V}n_{r'v} \le \gamma_{rir'i'} \le \sum_{v\in V}n_{rv} , \;\;\;\;\forall\; (r,i,r',i')\in ST \label{eq39}
\end{align}
Constraints \eqref{eq34}-\eqref{eq39} present the decisions related to time scheduling where constraints \eqref{eq35} is obtained after linearizing the term $ x_{rv}/u_{ri} $ in  constraints \eqref{eq10}. Replacing $ n_{r} $, the number of vessels on each liner service route $ r $ with the new variable $ n_{rv} $ (number of vessels of type $ v $ on each route $ r $), constraints \eqref{eq14} are transformed into the new set of constraints \eqref{eq39}.
\begin{align}
	& \sum_{\alpha\in A_{ri}}\psi_{ri\alpha}=1,\;\;\;\;\forall r\in R,\; \forall i\in I_r \label{eq40}\\
	& \sum_{\alpha\in A_{ri}}\alpha\psi_{ri\alpha}=u_{ri},\;\;\;\;\forall r\in R,\; \forall i\in I_r \label{eq41}\\
	&\phi_{riv\alpha} \le \psi_{ri\alpha},\;\;\;\;\forall r\in R,\;\forall i\in I_r,\;\forall v\in V,\;\forall \alpha\in A_{ri} \label{eq42}\\
	&\phi_{riv\alpha} \le x_{rv},\;\;\;\;\forall r\in R,\;\forall i\in I_r,\;\forall v\in V,\;\forall \alpha\in A_{ri} \label{eq43}\\
	&\phi_{riv\alpha} \ge \psi_{ri\alpha}+x_{rv}-1,\;\;\;\;\forall r\in R,\;\forall i\in I_r,\;\forall v\in V,\;\forall \alpha\in A_{ri} \label{eq44}\\
	& w_{riv\alpha} \le Cap_v \phi_{riv\alpha},\;\;\;\;\forall r\in R,\;\forall i\in I_r,\;\forall v\in V,\;\forall \alpha\in A_{ri} \label{eq61}
\end{align}
Constraints \eqref{eq40}-\eqref{eq61} are the newly defined constraints to linearize the non-linear terms related to vessel speed $u_{ri}$ and the binary variable $x_{rv}$.
\begin{align}
	&z_{riv} \le M_2x_{rv},\;\;\;\;\forall r\in R,\;\forall i\in I_r,\;\forall v\in V \label{eq45}\\
	&z_{riv} + M_3(1-x_{rv}) \ge \sum_{p\in P}\big(z^{load}_{rip}+z^{disc}_{rip}\big) ,\;\;\;\;\forall r\in R,\;\forall i\in I_r,\;\forall v\in V \label{eq46}\\
	& 168n_{rv}+ M_4(1-x_{rv}) \ge \sum_{i\in I_{r}}\big(T_{p_{ri}v}z_{riv} + \sum_{\alpha\in A_{ri}}L_{ri} \frac{\psi_{ri\alpha}}{\alpha}\big), \;\;\;\;\forall r\in R \;\forall v\in V \hspace{2cm} \label{eq47}
\end{align}
Constraints \eqref{eq45} and \eqref{eq46} define $ z_{riv} $, since these constraints ensure $ z_{riv} $ has feasible values if $ x_{rv} = 1 $.  Constraints \eqref{eq47} are the new constraints that imply ships have to follow weekly service frequency.
\begin{align}
	& \sum_{h\in H} \lambda_{rir'i'h}=1, \;\;\;\; \forall\; (r,i,r',i')\in ST \label{eq48}\\
	& \sum_{h\in H} h\lambda_{rir'i'h}=\theta_{rir'i'}, \;\;\;\; \forall\; (r,i,r',i')\in ST \label{eq49}\\
	& \delta_{prir'i'h} \ge \frac{1}{2} h\Big[\sum_{r\in R_p}\sum_{i\in I_{rp}}\sum_{\substack{o\in P\\o \ne p}}\big(z^{load}_{rio}+z^{disc}_{rio}\big)-\sum\limits_{d\in P}d_{pd}-\sum_{o\in P}d_{op}\Big]+ \notag\\
 &M_5(k_{prir'i'}+\lambda_{rir'i'h}-2),\;\;\;\;\;\;\;\;\;\;\;\;
	\forall p\in P, \;\forall\; (r,i,r',i')\in ST,\; \forall h\in H \label{eq50}
\end{align}
Equations \eqref{eq48}-\eqref{eq50} are the linearization constraints for the non-linear transshipment holding cost function in the objective function \eqref{eq1}.
\begin{align}
	&f_{rio}=0, \;\;\;\; \forall r\in R,\; \forall i\in I_r,\;o=p_{r,i+1} \label{eq51}\\
	&z^{disc}_{rio}=0, \;\;\;\; \forall r\in R,\; \forall i\in I_r,\;o=p_{ri} \label{eq52}\\
	&x_{rv},\;\phi_{riv\alpha},\;\psi_{ri\alpha}\in \{0,1\}, \;\;\;\;\forall r\in R,\;\forall i\in I_r,\;\forall v\in V,\;\forall \alpha\in A_{ri}, \label{eq53}\\
	& w_{riv\alpha} \ge 0, \;\;\;\;\forall r\in R,\;\forall i\in I_r,\;\forall v\in V,\;\forall \alpha\in A_{ri}, \label{eq62}\\
 & \lambda_{rir'i'h}\in \{0,1\}, \;\;\;\; \forall (r,i,r',i')\in ST,\;\forall h\in H \label{eq54}\\
	&n_{rv}\in \mathbb{Z}^+, \;\;\;\; \forall r\in R\label{eq55}\\
	&z^{disc}_{rio}\ge 0, z^{load}_{rio}\ge 0\;\;\;\; \forall r\in R,\; \forall i\in I_r,\; o\in P\label{eq56}\\
 	&z_{riv}\ge 0, \;\;\;\; \forall r\in R,\;\forall i\in I_r,\;\forall v\in V \label{eq57}
\end{align}
\begin{align}
	&f_{rio}\ge 0, \;\;\;\; \forall r\in R,\; \forall i\in I_r,\; o\in P \label{eq58}\\
	&\gamma_{rir'i'}\in \mathbb{Z},\;\theta_{rir'i'}\in \mathbb{Z}^+, \;\;\;\; \forall (r,i,r',i')\in ST \label{eq59}\\
	& t_{ri}, \; u_{ri}\in \mathbb{Z}^+, \;\;\;\; \forall r\in R,\; \forall i\in I_r \label{eq60}
\end{align}
Constraints \eqref{eq51} ensure that containerized cargos should not return to their origin port $ o$. Constraints \eqref{eq52} ensure that containerized cargos should not be discharged at their origin port $ o $. Finally, constraints \eqref{eq53}-\eqref{eq60} describe the integer attributes and non-negativity of the decision variables.
\par
Recently, an in-depth literature survey by \cite{christiansen2019liner} has shown that due to high complexity large-scale liner shipping problems are difficult to solve using exact methods and can lean on heuristics. Several popular conventional techniques such as goal programming, epsilon-constraint method,  and weighted sum method can be used to optimize a multiobjective problem. Nonetheless, the above-mentioned methods convert the multiobjective problem into single objective and then solve the single objective problem.
The above linear mathemtical model can be solved in the commercial softwares (e.g. CPLEX and LINGO) which are based on classical and conventional methods. But linearizing the model drastically raises the number of variables and constraints. Hence, solving linearized problem in commercial softwares has limitations such as, a large-scale problem will be very difficult to solve due to high computational complexity. In the original non-linear model, the sailing speed is considered as a continuous variable whereas in the linearized formulation it needs to be taken as a discrete variable. A finite discrete set, $ A_{ri} $ containing discrete speed levels has taken into account. We also have considered a discrete set of hours in a week ($ H $). 
\par
To overcome these limitations, Multi-objective Evolutionary algorithms (MOEAs) have a great potential in the domain of optimization to attain better quality solution within legitimate computational time \citep{tseng2008hybrid, zhang2023decomposition}.
Moreover, most of the classical techniques like goal programming, epsilon-constraint method,  and weighted sum method can yield single optimal solution in a single simulation run. MOEAs are able to provide multiple trade-off solutions along the Pareto frontier with modified selection schemes \cite{deb2001multi}. If a multi-objective optimization problem contains more than one optimal solutions, Evolutionary algorithms (EAs) can acquire the multiple pareto-optimal solutions in final iteration. Hence unlike conventional techniques, MOEAs find a set of trading off solutions in each generation.
In the domain of multi-objective optimization, MOEAs have attracted much attention of the researchers due to the following properties: (1) Finding multiple solutions forming optimal pareto front in only one simulation run, (2) Population based solution approach, (3) Elementary implementation and (4) Obtaining pareto-optimal solutions within reasonable computational time for large parameter search spaces and large-scale problem instances. 
\section{Solution Approach}\label{Solution approach}
	The multi-objective model with several decision variables and constraints becomes more complicated to solve when the number of routes, ports and ships increases. To avoid the computational complexity and to attain a near-optimal solution in such circumstances, the problem requires multi-objective EAs that can optimize multiple objectives simultaneously. Hence, two MOEAs: Nondominated Sorting Genetic Algorithm II (NSGA-II) and Online Clustering-based Evolutionary Algorithm (OCEA) are adopted to address the proposed problem.
	\subsection{Initial Solution}
	The MOEAs, NSGA-II and OCEA require an initial population or solution to start their search procedures. Initial population consists of all the decision variables which satisfy the corresponding constraints associated with the model. The problem contains different decision variables representing ship speed, optimal fleet, ship scheduling, and cargo handling. Constraints associated with the variables are considered and the initial solutions are generated correspondingly depending upon the constraints. The values of loading and unloading variables $ z^{load}_{rio} $ and  $ z^{load}_{rio} $ are computed by satisfying the demand constraint \eqref{eq8} along with constraints \eqref{eq16} \& \eqref{eq17}. Next, the obtained values are checked to satisfy the non-negativity constraints \eqref{eq20} \& \eqref{eq21} and discarded if infeasible. The variables representing flow of containers in each leg of the routes ($ f_{rio} $) are generated using the capacity constraint  \eqref{eq6} and the flow conservation constraint \eqref{eq7}. The fleet deployment variable, $ n_r $ is selected within a range given by \eqref{eq3} and ship speed, $ u_{ri} $ is chosen from the range given by constraint \eqref{eq5}. Similarly, the binary variable, $ x_{rv} $ is estimated such that it satisfies constraint \eqref{eq4}. Since the first vessel on each route serves the $ 1^{st} $ port of call within the first seven days, $ t_{r1} $ is generated by fulfilling constraint \eqref{eq9}. The arrival times of ships at each port of call except the first port ($  t_{r2},...,t_{ri},..$) are calculated using constraint \eqref{eq10}. Remaining variables $\gamma_{rir'i'}  $ and $ \theta_{rir'i'} $ are generated such that they satisfy constraints \eqref{eq13} and \eqref{eq14}. Hence, an initial solution is generated as discussed in the literature for various problems such as scheduling and routing. The initial solution thus achieved is provided to each algorithm so that the search procedures can be initiated.
	\subsection{Nondominated Sorting Genetic Algorithm II}
 
	A single-objective optimization problem can capture only one optimal solution, whereas a multi-objective optimization problem comes up with a set of non-dominated solutions i.e., the Pareto optimal front that depicts a trade-off between conflicting objectives. Nondominated Sorting Genetic Algorithm II (NSGA-II) is a commonly used MOEA that can achieve near optimal solutions handling more than one objectives together. This EA developed by \cite{NSGA-II} has an elite-preserving operator that gives the best population an opportunity to be transferred to the next generation.
	Among the MOEAs such as Nondominated Sorting Genetic Algorithm (NSGA), multi-objective genetic algorithm (MOGA), Strength Pareto Evolutionary Algorithm II (SPEA-II) etc., NSGA-II has become popular among the researchers due to their easy implementation and effectiveness \citep{rabbani2019stochastic}.
	NSGA-II outperformed Multi-Objective  Simulated Annealing (MOSA) and Multi-Objective Particle Swarm Optimization (MOPSO) in terms of efficiency, accuracy and speed of convergence for a three-objective closed-loop supply chain (CLSC) problem designed by \cite{babaveisi2018optimizing}. Using a modified NSGA-II, \cite{song2015multi} solved a stochastic multi-objective liner shipping problem to simultaneously optimize the service reliability, expected cost, and carbon emission under uncertain port time.
 Furthermore, \cite{de2017multiobjective} developed a bi-objective ship routing and service scheduling problem which enlightens port time window concept, ship draft restrictions, and sustainability. Two search heuristics, MOPSO and NSGA-II were applied to solve the complex engineering problem. Several researchers have proven its good solution quality and effectiveness by approaching combinatorial optimization problems through NSGA-II \citep{mogale2018grain}.

		\begin{figure}[H]
		\centering
				\frame{\includegraphics[width=0.5\linewidth]{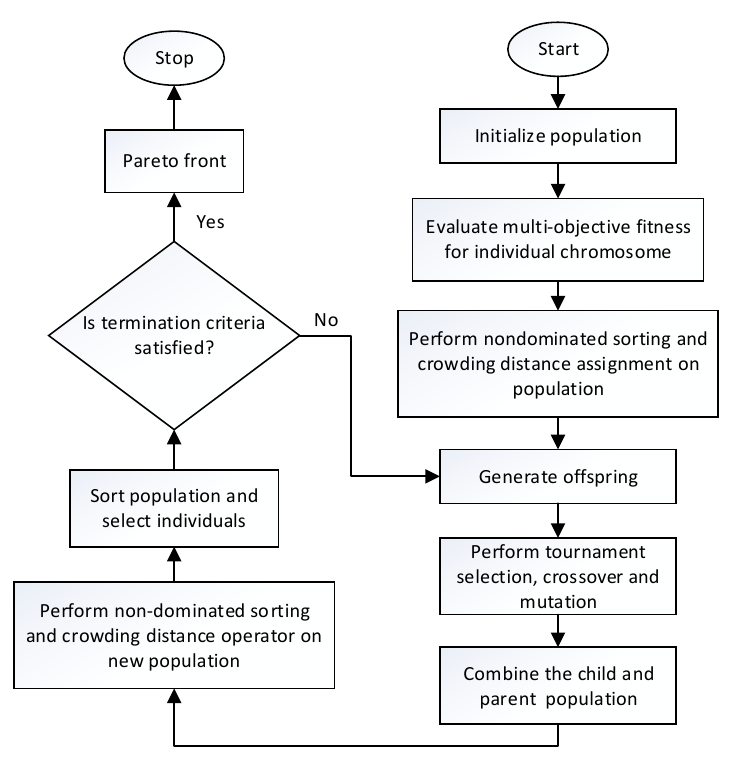}}
		\caption{Flowchart of NSGA-II}
		\label{fig:1}
		\end{figure}

	\subsection{Online Clustering-based Evolutionary Algorithm (OCEA)}
	A key element for developing an effective EA is to generate high-quality solutions. From the literature, it is evident that collecting information about the problem domain can significantly improve the ability to search \citep{zhang2008rm}.
	However, knowledge of an optimization problem's structure, in general, is unavailable. To gather information about the problem domain from the population or solutions, EA’s search pattern stimulates researchers to implement machine learning (ML) techniques. Hence, integrating machine learning methods with evolutionary algorithms is highly beneficial to obtain solutions of high quality. \cite{sun2018learning} proposed a learning-based MOEA with an online learning technique incorporated into the evolutionary search algorithm. While other learning-based MOEAs \citep{pelikan2006multiobjective,shim2013enhancing,li2014general} incur a high computational expense for learning, the online clustering-based evolutionary algorithm (OCEA) takes the series of solutions sequentially produced in the evolution as the training data.
	First, the solutions generated at each iteration are tested if they dominate the current population. If the new solutions do so, they will be introduced to the online learning method. Unlike batch-learning every non-dominated solution, regardless of its survival or elimination, is used only a single time for training.
	Therefore, we can consider the evolution data as a stream of data since each solution is created sequentially and not used for learning in next generations.
	\subsubsection{Algorithmic Framework}
	\begin{algorithm}[h]
		\caption{OCEA Framework}
		\label{alg1}
		\begin{algorithmic}[1]
			\Require{mating restriction probability ($ \alpha $), population size ($ K $), and maximum clusters considered ($ M_{max} $)}
			\Ensure{population of solutions ($ S $)}
			\State Initialize: population $ S = {s^1, . . . , s^K} $ and an external archive $ \mathcal{E} = S $.
			\State Form a cluster $ R^j $ with each $ s^j\in S $ where centroid $ c^j = s^j $ and
			counter $ h^j = 1 $, set $ R \leftarrow {R^j} $ and $ \tilde{R} = R $.
			\While {stopping criterion not satisfied}
			\State Let $ N \leftarrow |R| $.
			\State Choose a solution randomly from a $ R^j $, $ 1 \le j\le N $; create a global mating pool $ G $ from the selected solution.
			\For {$ j = 1 $ to $ K $}
			\State For each $ s^j $ create a mating pool $ M^j $:
			\vspace{-0.3cm}
			\begin{align}
				&M^i \leftarrow
				\begin{cases}
					R^{kj}\setminus \{s^j\} & \text{if }r_1<\alpha \\
					G & \text{otherwise}
				\end{cases} \notag\\
				\vspace{-0.3cm}
				&\text{where $ r_1=rand()$ generates a random number in $ [0,1] $ and $R^{kj} $ presents that} \notag\\
                &\text{$ s^j $ locates in $ R^{kj} $} \notag
			\end{align}
			\vspace*{-0.69cm}
			\State Find $ t^j\leftarrow SOLGEN(M^j, s^j) $.
			\State Clustering and modify $ [\mathcal{E},\tilde{R}] \leftarrow ESOC(\mathcal{E}, t^j, \tilde{R}) $
			\EndFor
			\State Put $ S\leftarrow \mathcal{E} $ and return the results of clustering ($ \tilde{R} $) of $ \mathcal{E} $ to $ S $
			\EndWhile
		\end{algorithmic}
	\end{algorithm}
	Algorithm \ref{alg1} depicts the pseudocode of OCEA. First, an initial population $ P $ is generated and set to be an initial external archive $ E $. The initial population is divided into clusters where every solution in the population is assumed to be a cluster and the centroid of the cluster is the solution itself, i.e. $ c^j = s^j $. The counter $ h^j $ of each cluster $ R^j $ is 1. Then total number of clusters is computed as $ N $ at each generation. A global mating pool $ G $ is created from a randomly selected solution of each $ R^j $. Also a mating pool $ M^j $ is constructed for each solution considering the mating restriction probability $ \alpha $, as shown in line 7. Next, a new solution $ t^j $ is evolved from the current solution $ s^j $ through the SOLGEN operator. Using the new solution, the cluster and external archive are updated. After $ K $ iterations, a new population is created for the subsequent generation using the updated cluster and external archive. Again, the same steps will be repeated until the termination criteria is satisfied.
	\subsubsection{New Solution Generation}
	\begin{algorithm}[h]
		\caption{Solution Generation Operator (SOLGEN)}
		\label{alg2}
		\begin{algorithmic}[1]
			\Require{an individual $ s $ and mating pool of the individual $ M $}
			\Ensure{a trial solution $ t $}	
			\State Select two different parent solutions $s^1 $ and $ s^2 $ from $ M $ randomly.	
			\State Create $ t' = (t_1',..., t_m')^T $ as following: for $ 1\le j\le m $,
			\vspace{-0.1cm}
			\begin{align}
				&t_j'=
				\begin{cases}
					s_j+C\times(s_j^1-s_j^2),& \text{if } r_2\le CR\\
					s_j,& \text{if } otherwise.
				\end{cases}\notag\\
				&\text{where $ r_2=rand()$ generates a random number in $ [0,1] $ and $CR$ is the }\notag\\
                &\text{crossover ratio.}\notag
			\end{align}
			\State Repair $ t': $ for $ 1\le j\le m $,
			\vspace{-0.1cm}
			\begin{align}
				t_j''=
				\begin{cases}
					l_j,& \text{if } t_j'<l_j\\
					u_j,& \text{if } t_j'>u_j\\
					t_j',& \text{if otherwise.}
				\end{cases}\notag
				\text{ where $ s_j \in [l_j, u_j] $.}\notag
			\end{align}
			\State Mutate $ t'': $ for $  1\le j\le m $,
			\vspace{-0.1cm}
			\begin{align}
				& 	t_j=
				\begin{cases}
					t_j''+\theta_j\times(u_j-l_j),& \text{if } r_3\le p\\
					t_j'',& \text{if otherwise.}
				\end{cases}\notag \\
				&\text{where $r_3 = rand()$ generates a random number in $[0,1]$, $p$ is the mutation }\notag\\
                &\text{probability and}\notag\\
				& \theta_j=
				\begin{cases}
					\bigg[2d+(1-2d)\Big(\frac{u_j-t_j''}{u_j-l_j}\Big)^{\rho+1}\bigg]^\frac{1}{\rho+1}-1,& \text{if } d<0.5\notag\\
					1-\bigg[2-2d+(2d-1)\Big(\frac{t_j''-l_j}{u_j-l_j}\Big)^{\rho+1}\bigg]^\frac{1}{\rho+1},& \text{if otherwise}.
				\end{cases} \notag\\
				&\text{where $d = rand()$ generates a random number in $[0,1]$ and $\rho$ is the }\notag\\
    			&\text{distribution index of mutation.}\notag
			\end{align}
			\State If necessary, repair $ t'' \leftarrow t $
		\end{algorithmic}
	\end{algorithm}
	For generating new solutions, OCEA applies two operators: polynomial mutation (PM) \citep{deb2001multi} and differential evolution (DE) \citep{price2006differential}. Algorithm \ref{alg2} represents the algorithm of new solution generation (SOLGEN()). This operator uses the present solution $ s $ and its mating pool $ M $ and results an offspring $ t $. A trial solution $ t' $ is produced from two parent solutions randomly chosen from mating pool using Differential Evolution (DE) where $ C $ is the scaling factor. Next, to ensure that the solution values lies within search range, a repair method is applied. Again, Polynomial Mutation (PM) is employed to create a new solution after mutating the repaired solution $ t'' $. If necessary, the solution is again repaired to be within boundary and then the final solution is sent back to the main algorithm.
 	At any generation, if the new solution is not dominated by the existing solutions of external archive, it is used for learning. Algorithm \ref{alg3} summarizes the procedure of ESOC(). First, $ \mathcal{E}\cup \{t\} $ is partitioned into $ I $ nondominated fronts $ {\mathcal{F}_1,...,\mathcal{F}_I} $ using elite nondomination sorting technique proposed by \cite{NSGA-II}. $ \mathcal{F}_1 $ is the best front and $ \mathcal{F}_I $ is the worst one. If the solutions in $ \mathcal{E}\cup \{t\} $ form more than one pareto front i.e. $ I > 1 $, the solution $ s^* $ in $ \mathcal{F}_I $ having highest $ d(s, \mathcal{E}\cup \{t\}) $ value is dropped. $ d(s,\mathcal{E}\cup \{t\}) $ is the number of solutions in $ \mathcal{E}\cup \{t\} $ that dominates $ s $. Otherwise, if $ I = 1 $, the solution $ s^* $ with minimum value of $ \Delta_\varphi(s,\mathcal{F}_1) $ is removed. $ \Delta_\varphi(s,\mathcal{F}) $ depicts the exclusive contribution of an individual $ s $ to the hypervolume measure or HV metric of its actual front $ \mathcal{F} $ \citep{beume2007sms} and is measured as $ \Delta_\varphi(s,\mathcal{F}) = \mathcal{H}(\mathcal{F})-\mathcal{H}(\mathcal{F}\setminus s)$ where $ \mathcal{H}(\mathcal{F}) $ denote the Lebesgue measure to the set $ \cup_{t\in \mathcal{F}}\{t'|t<t'<t_{ref}\} $ w.r.t a reference point $ t_{ref} $ and said to be the hypervolume (HV) metric. If $ s^*\ne t $ i.e., $ t $ is not discarded in the selection scheme, the online clustering mechanism is executed. First, $ s^* $ is discarded from its cluster $ R^* $, and the uodated values of centroid and counter of $ R^* $ are restored. Then a new cluster is formed with the centroid as $t $. If the total number of clusters in $ \mathcal{E} $ exceed $ M_{max} $, two clusters closest to each other are combined to finish the clustering.
	\subsubsection{Updating on Population and Clusters}
	\begin{algorithm}[h]
		\caption{Updating Procedure (ESOC)}
		\label{alg3}
		\begin{algorithmic}[1]
			\Require{external archive $\mathcal{E}$, allowed maximum number of clusters $
				N_{max}$, counters $h^n$ and centroids $c^n$ of present clusters $R^1,...,R^N$, $1\le n\le N$, and a new solution $t$.}
			\Ensure{external archive $\mathcal{E}$; clustering information of $\mathcal{E}$.}
			\State Implement the elitist non-dominated sorting technique on $\mathcal{E}\cup\{t\}$ acheive $I$ solution fronts ${\mathcal{F}_1,...,\mathcal{F}_I}$.
			\If {$I>1$}
			\State Find $s^*\leftarrow \displaystyle\arg\max_{s\in \mathcal{F}_I}d(s,\mathcal{E}\cup\{t\})$.
			\Else
			\State Find $s^*\leftarrow \displaystyle\arg\min_{s\in \mathcal{E}\cup \{t\}}\Delta_\varphi(s,\mathcal{F}_1)$.
			\EndIf
			\If {$s^*\neq t$}
			\State If $s^*\in R^n,n\in \{1,...,N\}$, then remove $s^*$ from $R^n: R^n\leftarrow R^n\setminus\{s^*\}$.
			\If {$R^n=\phi$}
			\State Remove $R^n$, set $N\leftarrow N-1$.
			\Else
			\State Update $R^n:h^n\leftarrow h^n-1$, $c^n\leftarrow c^n-\displaystyle\frac{s^*-c^n}{h^n}$.
			\EndIf
			\State Remove the worst solution $\mathcal{E}\leftarrow \mathcal{E}\cup \{t\}\setminus {s^*}$.
			\State Set $N\leftarrow N+1$, Make a new cluster $R^N$, put $h^N\leftarrow 1$, $c^N\leftarrow t$.
			\If {$N > N_{max}$}
			\State Determine: $(\xi,\eta)\leftarrow \displaystyle\arg\min_{\xi,\eta,\xi\neq\eta}\|c^\xi-c^\eta\|$
			\State Merge: $c^\xi\leftarrow \displaystyle\frac{c^\xi h^\xi+c^\eta h^\eta}{h^\xi+h^\eta}$, $h^\xi\leftarrow h^\xi+h^\eta$.
			\EndIf
			\Else
			\State Remove the worst solution in $\mathcal{E}\leftarrow \mathcal{E}\cup \{t\}\setminus \{s^*\}$.
			\EndIf
		\end{algorithmic}
	\end{algorithm}
		\section{Computational Study}\label{Computational Experiments}
	\subsection{Data description}
\begin{figure}[h]
		\centering
		\vspace*{-0.1cm}
		\begin{tikzpicture}
			\node (img)  {
				\fbox{\includegraphics[width=0.6\linewidth]{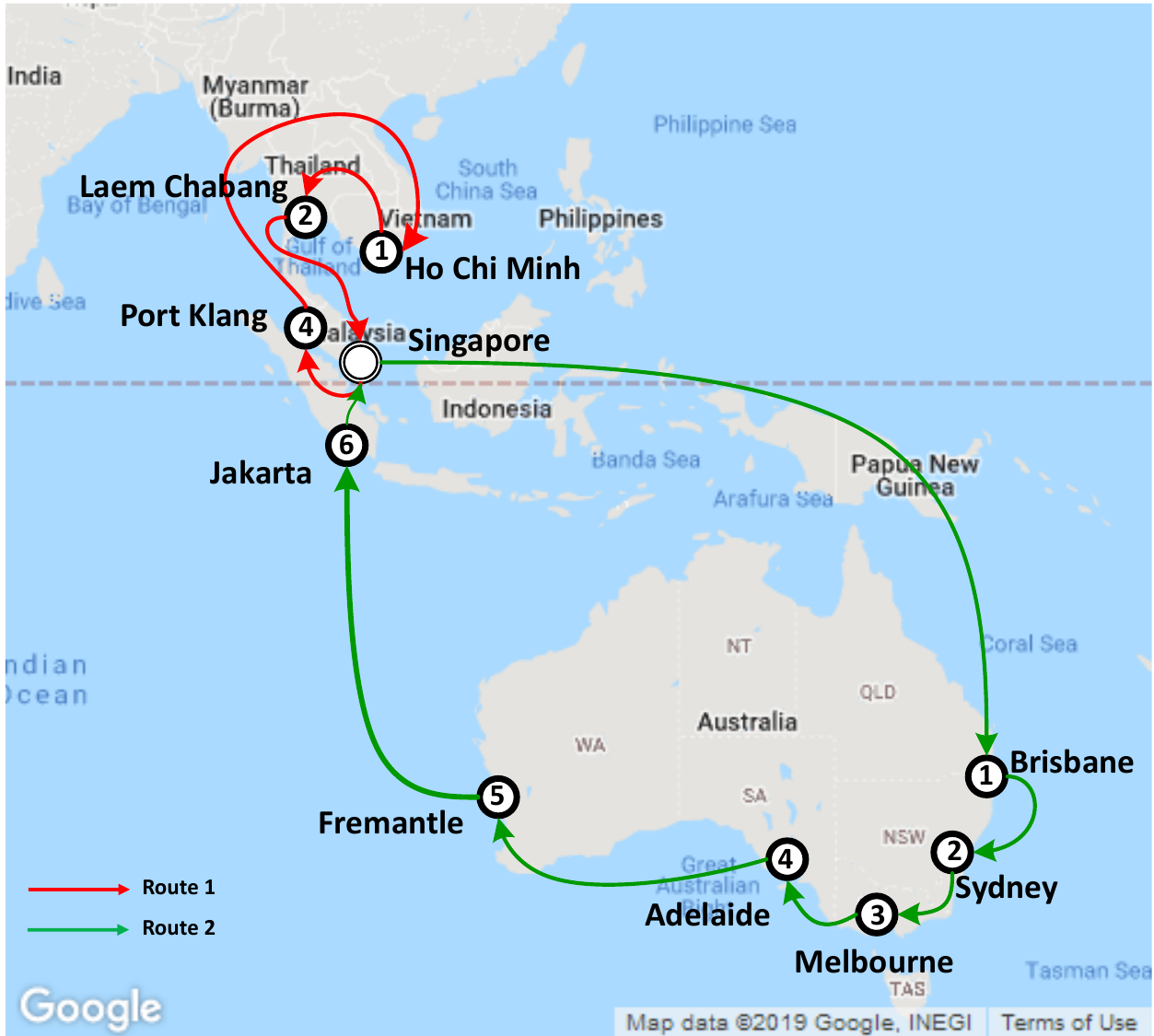}}};
		\end{tikzpicture}
		\caption{A liner service network containing two ship routes and a transshipment port}
		\vspace*{-0.3cm}
		\label{fig:2}
		\end{figure}

	To solve the proposed problem, six problem instances are simulated based on a real-life shipping network considering $ 6 $ ship routes and $ 24 $ ports. The routes contain four, eight, six, eight, six and eight ports of call respectively. Let us suppose that the ports have same characteristics apart from their position and containerized cargo demand. Table \ref{table1} presents the details about the ports of call and voyage length of each leg of the six routes. Two ship routes are shown in figure \ref{fig:2}. The loading cost $ C^{load}_p $ is assumed to be $ 150$ USD per TEU, unloading cost $ C^{disc}_p $ is $ 150$ USD per TEU, transshipped container handling cost $ C^{trans}_p $ is $ 150 $ USD per TEU \citep{meng2012liner}, and inventory holding cost for transshipped containers $ C^{hold}_p$ is $ 1.25 $ USD per TEU per hour. Five types of vessels are considered, as shown in table \ref{table2}. Fixed operating cost for each vesel, $ C^{opr}_v $ (USD), Berth occupancy charge for each vessel, $  C^{berth}_{pv} $, Fixed cost for calling each vessel at ports on each route,  $  C^{fix}_{4v} $, Capacity of each vessel, $Cap_v$ , Ship's lightweight with average weight of fresh water, ballast water, fuel, crew and provision, respectively (ton), and Container handling time at ports $ T_{pv}$ (hr/TEU) are given in table \ref{table2}. 	Minimum and maximum speed of ships ($u^{min}_{ri}$ and $u^{max}_{ri}$) on the $ i^\mathrm{th} $ leg on vessel route $ r $ are presumed to be 14 and 24 knots respectively \citep{xia2015joint}. Unit fuel cost, $c^{F}$ is 500 USD/ton, $ CO_2 $ emission factor in sea, $E_{sea}$ is 3.082 tons of $ CO_2 $/ton of fuel, $ CO_2 $ emission factor in ports $E_{port}$ is 0.01729 tons of $ CO_2 $/TEU handled and $ CO_2 $ emission cost $c^{E}$ is 32 USD/ton of $ CO_2 $ \citep{dulebenets2018greenijtst}.
 
	\begin{table}[h]
		\renewcommand*{\arraystretch}{1.4}
		\caption{Ports of call and length of each leg of the shipping routes}
		\fontsize{8pt}{8}\selectfont
		\centering
		\begin{tabular}{l l}
			\hline
			Sl no.  & Ports of call (length)  \\\hline
			1 & Ho Chi Minh (589) $ \rightarrow $ Laem Chabang (755) $ \rightarrow $ Singapore (187) $ \rightarrow $ Port Klang (830) $ \rightarrow $\\
            &Ho Chi Minh\\
			2 & Brisbane (419) $ \rightarrow $ Sydney (512) $ \rightarrow $ Melbourne (470) $ \rightarrow $ Adelaide (1325) $ \rightarrow $ Fremantle (1733) $ \rightarrow $\\
			& Jakarta (483) $ \rightarrow $ Singapore (3649) $ \rightarrow $ Brisbane.\\
			3 & Yokohama (15) $ \rightarrow $ Tokyo (177) $ \rightarrow $ Nagoya (201) $ \rightarrow $ Kobe (734) $ \rightarrow $ Shanghai (745) $ \rightarrow $\\
			& Hong Kong (1568) $ \rightarrow $ Yokohama\\
			4 & Dalian (187) $ \rightarrow $ Xingang (379) $ \rightarrow $ Qingdao (303) $ \rightarrow $ Xiamen (93) - Ningbo (93) $ \rightarrow $ Shanghai\\
			&(383) $ \rightarrow $ Kwangyang (72) $ \rightarrow $ Busan (487) $ \rightarrow $ Dalian\\
			5 & Ho Chi Minh (589) $ \rightarrow $ Laem Chabang (755)$ \rightarrow $ Singapore (187)$ \rightarrow $ Port Klang (830) $ \rightarrow $ Qingdao\\
			& (345) $ \rightarrow $ Shanghai (876) $ \rightarrow $ Ho Chi Minh\\
			6 & Brisbane (419)$ \rightarrow $ Sydney (512) $ \rightarrow $ Melbourne (470)$ \rightarrow $ Adelaide (1325)$ \rightarrow $ Fremantle (1733)$ \rightarrow $\\
			& Jakarta (483) $ \rightarrow $ Singapore (3649)$ \rightarrow $ Colombo (1287) $ \rightarrow $ Brisbane\\\hline
		\end{tabular}\label{table1}
	\end{table}
	\begin{table}[h]
		\renewcommand*{\arraystretch}{1.4}
		\caption{Characteristics of the liner vessels}
		\fontsize{8pt}{8}\selectfont
		\centering
		\begin{tabular}{l ccccc}
			\hline
			\multirow{2}{*}{Parameter Values} & \multicolumn{5}{c}{Vessel type}  \\\cline{2-6}
			& 1 & 2 & 3 & 4 & 5    \\\hline
			Fixed operating cost, $ C^{opr}_v $ (USD) &37,485 & 51,923& 76,923& 115,384 & 173,076 \\
			Berth occupancy charge, $  C^{berth}_{pv} $(USD/h) & 500 & 1000& 1666 & 3333 & 5000\\
			Fixed cost for calling at ports on route 1, $  C^{fix}_{1v} $ (USD) & 154791 & 191900 &240500 &256600& 276100 \\
			Fixed cost for calling at ports on route 2, $  C^{fix}_{2v} $ (USD) & 533980 & 689651 &788300& 854600 &929100 \\
			Fixed cost for calling at ports on route 3, $  C^{fix}_{3v} $ (USD) & 226198 & 280542 &342760 &384500& 404000 \\
			Fixed cost for calling at ports on route 4, $  C^{fix}_{4v} $ (USD) & 148807 & 187600 & 220850 & 259800 &279700 \\
			Fixed cost for calling at ports on route 5, $  C^{fix}_{5v} $ (USD) & 197892 & 235340 & 292760 & 304500 & 324000 \\
			Fixed cost for calling at ports on route 6, $  C^{fix}_{6v} $ (USD) & 594070 & 730527 & 840582 & 929753 & 989650 \\
			Capacity, $Cap_v$ (TEUs) &2400 & 4800 & 8400 & 11000 & 15000 \\
			Ships' lightweight + average weight of fuel, water, & 21832 & 36898 & 54,753	& 66,204 & 79,612\\
            provision, crew (ton) \\
			Fixed time when calling at a port (hr) &4 &4 &4 & 4 & 4\\ 
			Container handling time $ T_{pv}$ (hr/TEU) & 0.025 & 0.012 & 0.011 & 0.008 & 0.007 \\
			Maximum number of ships $n^{max}_r$ & 15 & 15 & 15 &15 & 15\\\hline
		\end{tabular}\label{table2}
	\end{table}
	\subsection{Results and discussion}
\begin{figure}[h]
		\centering
		\begin{minipage}[b]{0.441\linewidth}
			\centering
			\fbox{
				\includegraphics[width=1.04\linewidth]{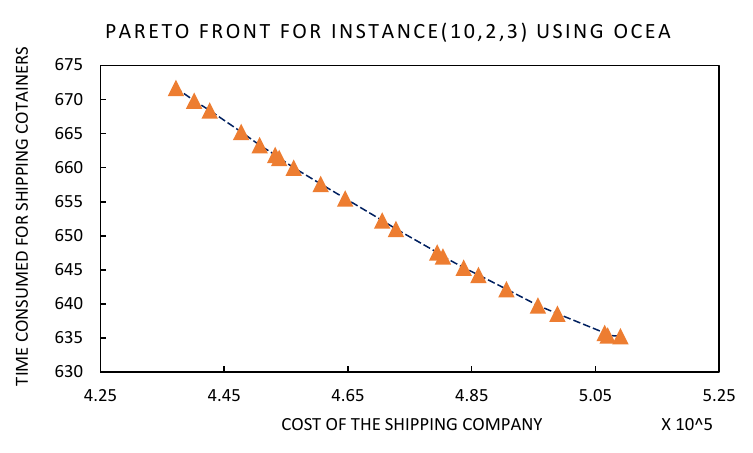}}
			\caption{Pareto front for instance (10,2,3) using OCEA}
			\label{fig:3}
		\end{minipage}
		\qquad
		\begin{minipage}[b]{0.45\linewidth}
			\centering
			\fbox{\includegraphics[width=1.04\linewidth]{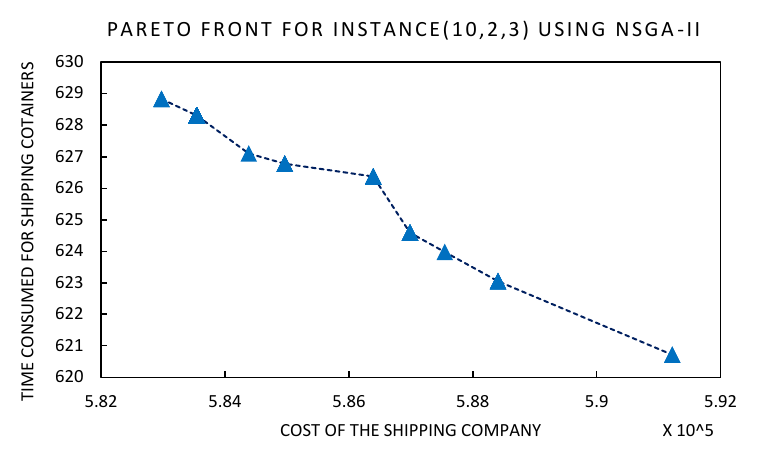}}
			\caption{Pareto front for instance(10,2,3) using NSGA-II}
			\label{fig:4}
		\end{minipage}
	\end{figure}

\begin{figure}[h]
		\centering
		\begin{minipage}[b]{0.44\linewidth}
			\centering
			\fbox{
				\includegraphics[width=1.04\linewidth]{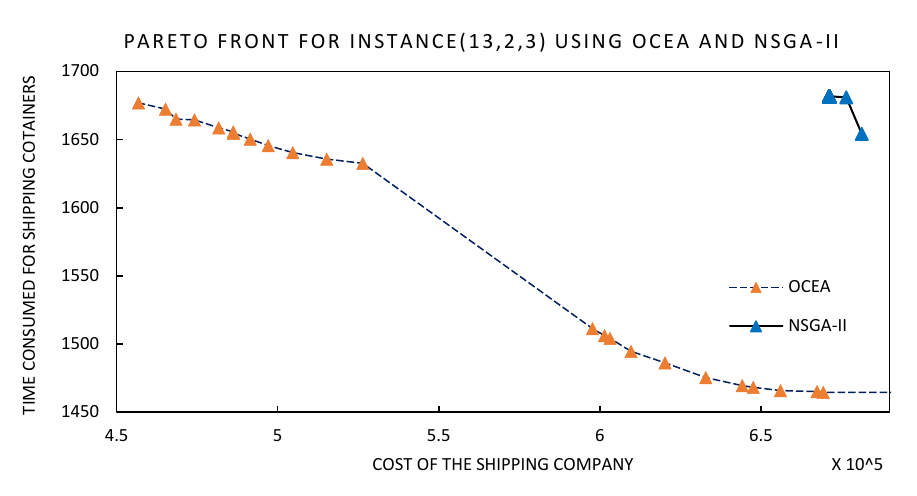}}
			\caption{Pareto front for instance (13,2,3) using OCEA and NSGA-II}
			\label{fig:5}
		\end{minipage}
		\qquad
		\begin{minipage}[b]{0.43\linewidth}
			\centering
			\fbox{\includegraphics[width=1.04\linewidth]{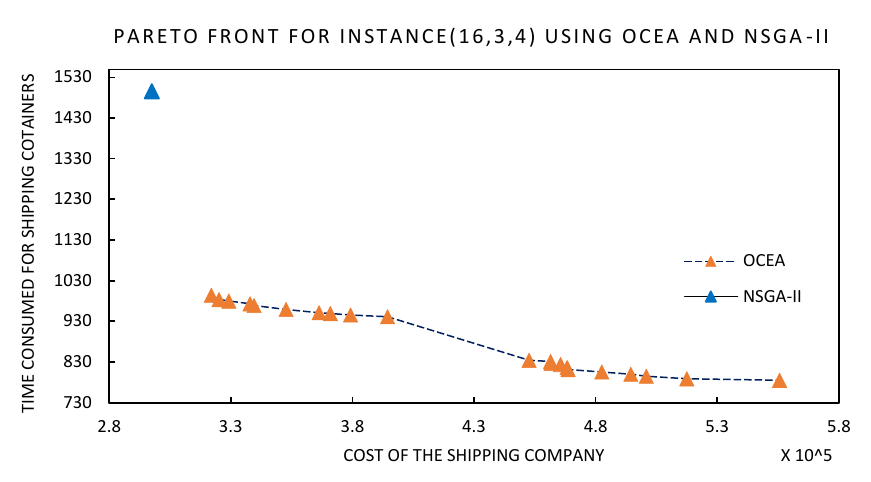}};
			\caption{Pareto front for instance (16,3,4) using OCEA and NSGA-II}
			\label{fig:6}
		\end{minipage}
	\end{figure}

	Numerical experiments based on the test problems are performed to verify and solve the multi-objective model proposed in this study. The six problem instances are generated making variation in the total number of ports, number of routes and number of different type of vessels. The computational experiments are performed on a computer with Intel Core i5, 3.20GHz processor and 4.00 GB installed physical memory (RAM). To validate and solve the test problems, Online Clustering-based Evolutionary Algorithm (OCEA) and NSGA-II: two EAs are exploited in this study and coded in MATLAB (R2018b) software. The best parameter settings of the algorithms are obtained after a set of test runs. For OCEA, the population size is set to 100; mutation probability is $p = 1/m$; distribution index of mutation $\rho = 20$; maximum number of clusters $ N_{max} = 3$; $\beta = 0.6$; scaling factor $C=1$; and crossover ratio $CR= 0.5$. For NSGA-II, the population size = 100; mutation probability is $1/m$; distribution index of mutation $= 100$; crossover ratio $= 0.5$; and distribution index of crossover $= 20$. Table \ref{table3} discusses the complexities of six problem instances related to the number of variables and number of constraints. It also presents the computational efficiency of OCEA and NSGA-II to solve all the problem instances.
	\begin{table}[h]
		\renewcommand*{\arraystretch}{1.4}
		\caption{Different problem instances and elapsed time to solve the test problems}
		\fontsize{8pt}{8}\selectfont
		\centering
		\begin{tabular}{cccccc}
			\hline
			Sl & Problem instances  &  No. of & No. of  & Computational time & Computational time \\
			no.	  & (ports, routes,  &  variables &
			constraints  & elapsed in & elapsed in \\
            	  & vessels)  &   &
			  & NSGA-II (min) & OCEA (min) \\\hline
			1 & (10,2,3) & 604 & 1049 & 12 & 171 \\
			2 & (13,2,3) & 976 & 1797 & 13 & 664\\
			3 & (16,3,4) & 1444 & 2726 & 32 &  1248\\
			4 & (18,3,4) & 1804 & 3435 & 30 & 1552 \\
			5 & (24,4,5) & 3128 & 6334 & - & 1810 \\
			6 & (27,4,5) & 3920 & 7975 & - & 1890 \\
			\hline
		\end{tabular}\label{table3}
	\end{table}
	\begin{table}[h]
		\renewcommand*{\arraystretch}{1.4}
		\caption{Cost incurred and total voyage time for each of the problem instances}
		\fontsize{8pt}{8}\selectfont
		\centering
		\begin{tabular}{c c ccccc}
			\hline
			Serial & Problem instances  & \multicolumn{2}{c} {Shipping cost incurred using} & & \multicolumn{2}{c}{Shipping time using}\\\cline{3-4}\cline{6-7}
			no.	  & (ports, routes, vessels)  & NSGA–II (USD) & OCEA (USD)  & & NSGA–II (hour)  &  OCEA (hour) \\\hline
			1 & (10,2,3) &$ 5.829\times 10^5 $ & $ 5.09\times10^5 $ & & $ 6.288\times10^2 $ & $ 6.352\times10^2 $\\
			2 & (13,2,3) & $ 6.713\times 10^5 $ & $ 4.568\times 10^5 $ & & $ 1.681\times 10^3 $ & $ 1.677\times 10^3 $\\
			3 & (16,3,4) & $ 2.974\times 10^6 $ & $ 3.219\times 10^6 $ & &$ 1.495\times 10^3 $ &$ 0.994\times 10^3 $\\
			4 & (18,3,4) & $ 7.472\times 10^6 $ & $ 4.213\times 10^6 $ & &$ 1.060\times 10^3 $  &$ 0.982\times 10^3 $ \\
			5 & (24,4,5) & & $ 3.734\times 10^5 $ & & & $ 1.417\times 10^3 $ \\
			6 & (27,4,5) & & $ 2.915\times 10^6 $ & & & $ 0.665\times 10^3 $\\
			\hline
		\end{tabular}\label{table4}
	\end{table}
	\begin{table}[h]
		\renewcommand*{\arraystretch}{1.4}
		\caption{Pareto-optimal solutions for the instances $(13,2,3)$ attained using OCEA and NSGA–II}
		\fontsize{8pt}{8}\selectfont
		\centering
		\begin{tabular}{c ccccc}
			\hline
			No. of Pareto & \multicolumn{2}{l}{NSGA–II results: Problem instance $(13,2,3)$}  & & \multicolumn{2}{l}{OCEA results: Problem instance $(13,2,3)$}\\\cline{2-3}\cline{5-6}
			solutions  & Shipping cost incurred & Shipping time & & Shipping cost incurred & Shipping time \\\hline
			1 & $ 6.71307\times 10^5 $	& $ 1.682\times 10^3 $ & & $ 4.56786\times 10^5 $ & $ 1.677\times 10^3 $ \\	
			2 & $ 6.76548\times 10^5 $  & $ 1.681\times 10^3 $  & & $ 4.68465\times 10^5 $ & $ 1.665\times 10^3 $\\	
			3 & $ 6.81345\times 10^5 $	& $ 1.654\times 10^3 $ & & $ 4.81753\times 10^5  $& $ 1.658\times 10^3 $ \\	
			4 &                         &                      & & $ 4.91485\times 10^5 $ & $ 1.650\times 10^3 $ \\
			5 &                     	&                     & & $ 5.15280\times 10^5 $ & $ 1.635\times 10^3 $ \\
			6 &                     	&                     & & $ 5.97778\times 10^5 $ & $ 1.511\times 10^3 $ \\\hline
		\end{tabular}\label{table5}
	\end{table}
	For each test problem, one solution from the non-dominated solution set of the pareto front ranked 1 is chosen and provided in table \ref{table4}. One optimal non-dominated solution in terms of total cost incurred and shipping time are determined using OCEA and NSGA-II, as given in the table \ref{table4}. For example, the total shipping costs for the test problem (13,2,3) obtained after applying the NSGA-II and online clustering-based EA are  $ 6.713\times 10^5 $ and $ 4.568\times 10^5 $ USD respectively. Similarly, total shipping times for the test problem (16,3,4) are $ 1.495\times 10^3 $ and $ 0.994\times 10^3 $ hours when solved using NSGA-II and OCEA respectively. In case of the problem instances (24,4,5) and (27,4,5), NSGA-II is unable to find the feasible solutions even after searching for large number of iterations whereas the optimal shipping time and cost incurred can be easily attained using OCEA. From the table \ref{table4}, it can be summarized that OCEA can find the pareto optimal solutions more efficiently than NSGA-II. Figures \ref{fig:3} and \ref{fig:4} graphically portray the solutions of the pareto optimal front for the test case (10,2,3) using OCEA and NSGA-II respectively. Next in problem instance (13,2,3), the trade-off between total cost and shipping time is presented in figure \ref{fig:5} for both the evolutionary algorithms. Similarly, table \ref{table5} depicts all the three pareto solutions obtained using NSGA-II for the problem instance (13,2,3) and only six solutions among all the pareto-optimal solutions obtained after applying the online clustering-based EA. It can be seen that from the figure \ref{fig:5} and table \ref{table5}, OCEA can capture optimal solutions more efficiently and also the number of pareto solutions obtained are much higher than NSGA-II.
	
\begin{figure}[h]
	\centering
	\begin{minipage}[b]{0.44\linewidth}
		\centering
		\fbox{\includegraphics[width=1.04\linewidth]{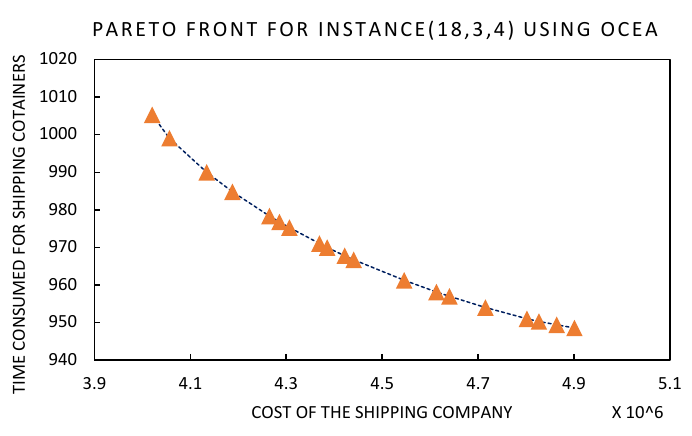}}
		\caption{Pareto front for instance $ (18,3,4) $ using OCEA}
		\label{fig:7}
	\end{minipage}
	\qquad
	\begin{minipage}[b]{0.438\linewidth}
		\centering
		\fbox{\includegraphics[width=1.04\linewidth]{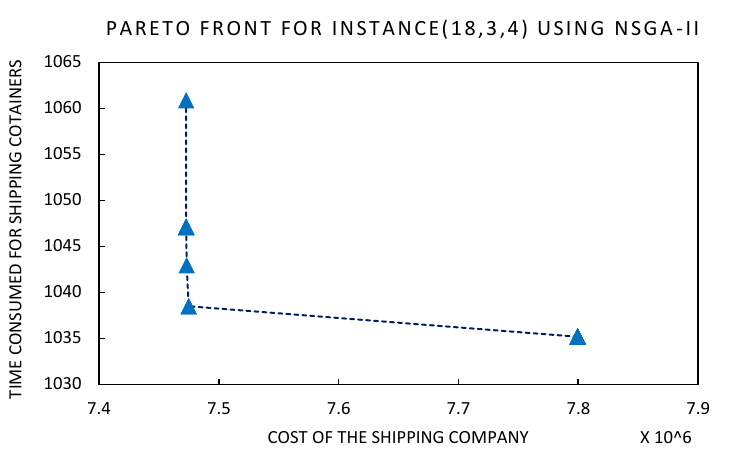}}
		\caption{Pareto front for instance $ (18,3,4) $ using NSGA-II}
		\label{fig:8}
	\end{minipage}
\end{figure}

\begin{figure}[H]
		\centering
		\begin{minipage}[b]{0.446\linewidth}
			\centering
			\fbox{
				\includegraphics[width=1.04\linewidth]{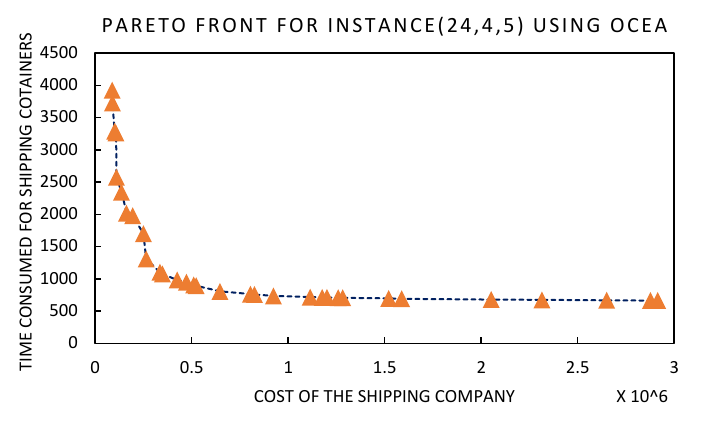}}
			\caption{Pareto front for instance $ (24,4,5) $ using OCEA}
			\label{fig:9}
		\end{minipage}
		\qquad
		\begin{minipage}[b]{0.45\linewidth}
			\centering
			\fbox{\includegraphics[width=1.04\linewidth]{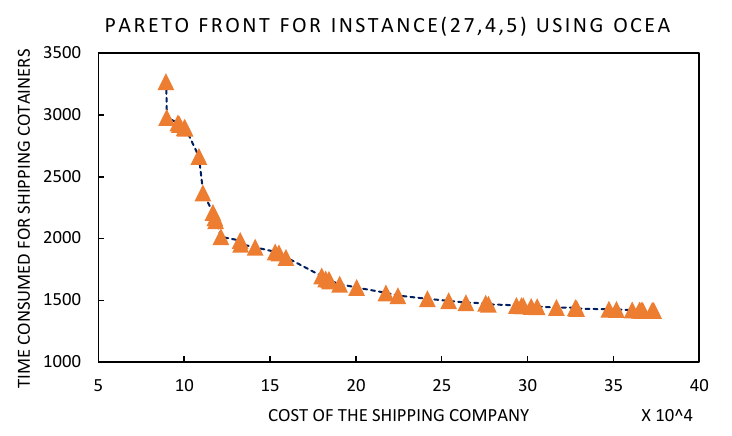}}
			\caption{Pareto front for instance $ (27,4,5) $ using OCEA}
			\label{fig:10}
		\end{minipage}
	\end{figure}

	Multiple pareto optimal non-dominated solutions hold good quality and diversity in decision variables and objective. For the problem instance $ (16,3,4) $, the trade-off points of the global pareto optimal front acheived with OCEA and NSGA-II are captured in figure \ref{fig:6}.
	These figures reveal that if a liner shipping company accelerates the vessel speed to consume less time, then the overall shipping cost will increase with a significant increase in carbon emissions (since carbon emission is considered in terms of a cost factor here). Hence, the near-optimal pareto front solutions of this multiobjective model will assist managers to manitain a balance between overall cost, carbon emission and voyage time by exploring different combinations of the solutions and objectives as per the demand. Furthermore, table \ref{table6} represents the pareto optimal front containing five non-dominated solutions 
 found using NSGA-II for the test case $ (18,3,4) $. Among all the solutions of the first pareto front acheived using OCEA for problem instance $ (18,3,4) $, twelve non-dominated solutions are considered and mentioned in table \ref{table6}, whereas figure \ref{fig:7} portrays all the solution points of the first pareto front while using the same algorithm. Figure \ref{fig:8} represents the first pareto front solution points for the same problem using NSGA-II. Finally, figures \ref{fig:9} and \ref{fig:10} present the large instances $(24,4,5)  $ and $ (27,4,5) $ respectively and their optimal pareto front containing all the non-dominated solutions attained using online clustering based EA. Since the EAs include diversity-preserving operator, the pareto-optimal set gives good quality solutions that are diverse in both the objective functions and hence, multiple possibilities with different solutions also help the shipping company manager to make an appropriate choice.
 
		\begin{table}[h]
		\renewcommand*{\arraystretch}{1.4}
		\caption{Pareto-optimal solutions for the instances $(18,3,4)$ attained using OCEA and NSGA–II}
		\fontsize{8pt}{8}\selectfont
		\centering
		\begin{tabular}{c ccccc}
			\hline
			No. of Pareto & \multicolumn{2}{l}{NSGA–II results: Problem instance $(18,3,4)$}  & & \multicolumn{2}{l}{OCEA results: Problem instance $(18,3,4)$}\\\cline{2-3}\cline{5-6}
			solutions	  & Shipping cost incurred & Shipping time & & Shipping cost incurred & Shipping time \\\hline
			1 & $ 7.47257\times 10^6 $	& $ 1.060\times 10^3 $ & & $ 4.02046\times 10^6 $ & $ 1.005\times 10^3 $\\	
			2 & $ 7.47284\times 10^6 $  & $ 1.047\times 10^3 $  & & $ 4.08085\times 10^6 $ &	 $0.995\times 10^3 $\\	
			3 & $ 7.47309\times 10^6 $	& $ 1.042\times 10^3 $ & & $ 4.13382\times 10^6  $&  $0.990\times 10^3 $\\	
			4 & $ 7.47488\times 10^6 $  & $ 1.038\times 10^3 $ & & $ 4.21364\times 10^6 $ & $0.982\times 10^3 $ \\
			5 & $ 7.79949\times 10^6 $ 	& $ 1.035\times 10^3 $& & $ 4.26485\times 10^6 $ & $0.978\times 10^3 $ \\
			6 &                     	&         & & $ 4.33607\times 10^6 $ & $ 0.973\times 10^3 $ \\ 	
			7 & 	                    &		  & & $ 4.40209\times 10^6 $ & $ 0.969\times 10^3 $ \\
			8 & 	                    & 		  & & $ 4.66023\times 10^6 $ & $ 0.956\times 10^3 $\\
			9 & 	                    & 		  & & $ 4.50591\times 10^6 $ & $ 0.963\times 10^3 $\\
			10 & 	                    & 		  & & $ 4.59663\times 10^6 $ & $ 0.958\times 10^3 $ \\
			11 & 	                    & 		  & & $ 4.78362\times 10^6 $ & $ 0.951\times 10^3 $ \\
			12 & 	                    & 		  & & $ 4.86412\times 10^6 $ & $ 0.949\times 10^3 $ \\\hline	
		\end{tabular}\label{table6}
	\end{table}
 
\vspace{-0.5cm}
	\subsection{Managerial implication}
	For high revenue generation, shipping companies intend to lower the total shipping time and as a consequence they increase the sailing speed of ships causing more fuel consumption which is approximately cubic in speed. Higher fuel consumption results in larger carbon emissions. There is an urgent need to maintain an equilibrium between overall cost including carbon emissions and shipping time. From this study it can be stated that considering carbon emission with several realistic factors such as transshipment cost and payload-dependent fuel consumption in the overall shipping cost is an integral part of the multi-objective liner shipping problem for capturing the sustainability issues as well as the real-life complexities of seaborne shipping. There has been a growing interest in the past few years to adopt sustainability management practices within the domain of supply chain (cite some papers from selected journal). Multi-objective optimization seems to be a viable approach to capture the conflicts between economic and environmental sustainability in the modeling of strategic and operational decisions in seaborne transportation. Overall, this mathematical model will help the decision manager to take several interconnected decisions associated with transshipping containers, deployment of fleet, sailing speed optimization, vessel scheduling and cargo flow so that sustainability of the maritime logistics enhances.
	
	\section{Conclusion}\label{Conclusion}
	For liner shipping companies, a mathematical multi-objective model is proposed here to simultaneously plan the suitable service schedule, number of vessels in a fleet, vessel speed on each leg, and the containerized cargo allocation for each pair of origin-destination. This multi-objective MINLP model presents a trade-off between economical and environmental aspects considering total shipping time and overall shipping cost including carbon emission and several operational costs as the two conflicting objectives. The major contributions of the work are threefold: (1) This is a holistic approach that can be applied to a generic liner shipping service with more than one routes with transshipment hubs.
 The objective functions include various operational costs, all of which can satisfy the major interests of a liner shipping manager. (2) Fuel consumption function considered here is speed and payload dependent where a limited number of studies have included payload in their model. (3) This model also captures the environmental aspect by integrating the carbon emission cost in the first objective. However, the cargo allocation decision is connected with the service schedule and fleet deployment decisions. Numerous researchers have considered all these decision variables in their work, but our study has simultaneously taken care of all major decisions. These decisions are highly connected with each other, which makes the proposed model non-linear and much complex. NSGA-II and OCEA: two intelligent and efficient evolutionary algorithms are applied to attain the near-optimal solution of the proposed problem. Furthermore, six problem instances of different sizes are solved using these algorithms to validate the proposed model. Managerial insights drawn from the model will lead the shipping companies to improve their operations management policies.
	\par
	Future studies can further consider uncertainties in port operations, such as uncertain wait time for ships because of congested ports and uncertain port stay time for handling containers in the liner shipping problem. Although this model incorporates the bunker price as a constant, in future it can be considered as stochastic. These uncertainties have a significant impact on scheduling ships' visit time at ports, and bunker consumption and therefore, the vessel speed on each sailing leg. Furthermore, uncertain container shipping demand also can be included in the model as an interesting future research.

 \section{Disclosure statement}
 The authors have no conflicts of interest.

 \section{Data availability statement}
 The data that support the findings of this study are available from the corresponding author, Jasashwi Mandal, upon reasonable request.

\bibliography{References_report}
\end{document}